\documentclass{amsart}

\usepackage{amssymb}
\usepackage{amsfonts,amsmath,amssymb,mathrsfs,verbatim}
\usepackage[dvips]{graphicx}
\usepackage{graphics}
\usepackage{psfrag}
\usepackage{hyperref}
\usepackage{marginnote}
\usepackage{xcolor}

\begin{document}

\renewcommand{\P}{\mathbb{P}}
\newcommand{\N}{\mathbb{N}}
\newcommand{\A}{\mathbb{A}}
\newcommand{\Z}{\mathbb{Z}}
\newcommand{\Q}{\mathbb{Q}}
\newcommand{\R}{\mathbb{R}}
\renewcommand{\d}{\underline{d}}
\newcommand{\lcm}{\mathrm{lcm}}
\newcommand{\mA}{\mathcal{A}}
\newcommand{\mB}{\mathcal{B}}
\newcommand{\mC}{\mathcal{C}}
\newcommand{\mH}{\mathcal{H}}
\newcommand{\mS}{\mathcal{S}}
\newcommand{\mP}{\mathcal{P}}
\newcommand{\mM}{\mathcal{M}}
\newcommand{\mK}{\mathcal{K}}
\newcommand{\mI}{\mathcal{I}}
\newcommand{\mU}{\mathcal{U}}
\newcommand{\mT}{\mathcal{T}}
\newcommand{\IA}{I_{\mA}}
\newcommand{\IB}{I_{\mB}}
\newcommand{\IS}{I_{\mS}}
\newcommand{\IM}{I_{\mM}}

\newcommand{\lS}{\leq_{\mS}}
\newcommand{\x}{\mathbf{x}}
\renewcommand{\t}{\mathbf{t}}
\newcommand{\kx}{k[\x]}
\newcommand{\kt}{k[\t]}
\newcommand{\mG}{\mathcal{G}}

\newtheorem{thm1}{Theorem}[section]
\newtheorem{lem1}[thm1]{Lemma}
\newtheorem{rem1}[thm1]{Remark}
\newtheorem{def1}[thm1]{Definition}
\newtheorem{cor1}[thm1]{Corollary}
\newtheorem{defn1}[thm1]{Definition}
\newtheorem{prop1}[thm1]{Proposition}
\newtheorem{ex1}[thm1]{Example}
\newtheorem{alg1}[thm1]{Algorithm}
\newtheorem{prob1}[thm1]{Open problem}
\newtheorem{quest1}[thm1]{Question}

\title{On robustness and related properties on toric ideals}
\author{Ignacio Garc\'ia-Marco}
\address{Ignacio Garc\'ia-Marco, Facultad de Ciencias, Universidad de La Laguna, 38200 La Laguna, Tenerife, Spain}
\email{iggarcia@ull.edu.es}
\author[Christos Tatakis]{Christos Tatakis$^{*}$}\thanks{* Corresponding author.}
\address{Christos Tatakis, Department of Mathematics, University of Ioannina, 45110 Ioannina, Ioannina, Greece}
\email{chtataki@uoi.gr}


\begin{abstract}
A toric ideal is called robust if its universal Gr\"obner basis is a minimal set of generators, and is called generalized robust if its universal Gr\"obner basis  equals  its universal Markov basis (the union of all its minimal sets of binomial generators). Robust and generalized robust toric ideals are both interesting from both a Commutative Algebra and an Algebraic Statistics perspective. However, only a few nontrivial examples of such ideals are known. In this work we study these properties for toric ideals of both graphs and numerical semigroups. 
For toric ideals of graphs, we characterize combinatorially the graphs giving rise to robust and to generalized robust toric ideals generated by quadratic binomials. As a byproduct, we obtain families of Koszul rings. For toric ideals of numerical semigroups, we determine that one of its initial ideals is a complete intersection if and only if the semigroup belongs to the so-called family of free numerical semigroups. Hence, we characterize all complete intersection numerical semigroups which are minimally generated by one of its Gr\"obner basis and, as a consequence, all the Betti numbers of the toric ideal and its corresponding initial ideal coincide. Moreover, also for numerical semigroups, we prove that the ideal is generalized robust if and only if the semigroup has a unique Betti element and that there are only trivial examples of robust ideals. We finish the paper with some open questions.  
\end{abstract}

\maketitle
\tableofcontents

\section{Introduction}

Let $A=\{\textbf{a}_1,\ldots,\textbf{a}_m\}\subseteq \mathbb{N}^n$
be a finite set of nonzero vectors  and
$\mathbb{N}A:=\{l_1\textbf{a}_1+\cdots+l_m\textbf{a}_m \ | \ l_i \in\mathbb{N}\}$ the 
corresponding affine monoid. We grade the
polynomial ring $\mathbb{K}[x_1,\ldots,x_m]$ over an arbitrary field $\mathbb{K}$ by the
semigroup $\mathbb{N}A$ setting $\deg_{A}(x_i)=\textbf{a}_i$ for
$i=1,\ldots,m$. For $\textbf{u}=(u_1,\ldots,u_m) \in \mathbb{N}^m$,
we define the $A$-\emph{degree} of the monomial $\textbf{x}^{\textbf{u}}:=x_1^{u_1} \cdots x_m^{u_m}$
to be $\deg_{A}(\textbf{x}^{\textbf{u}}):=u_1\textbf{a}_1+\cdots+u_m\textbf{a}_m
\in \mathbb{N}A,$
while we denote the usual degree $u_1+\cdots +u_m$ of 
$\textbf{x}^{\textbf{u}}$ by $\deg(\textbf{x}^{\textbf{u}})$. The
\emph{toric ideal} $I_{A}$ associated to $A$ is the ideal generated by all the binomials $\textbf{x}^{\textbf{u}}- \textbf{x}^{\textbf{v}}$
such that $\deg_{A}(\textbf{x}^{\textbf{u}})=\deg_{A}(\textbf{x}^{\textbf{v}})$. It is a prime ideal of height $m - r$, being $r$ the rank of the subgroup of $\Z^m$ spanned by $A$ (see, e.g., \cite{ST}).  Toric ideals have applications in several areas such as: algebraic statistics, biology, computer algebra,
computer aided geometric design, dynamical systems,
hypergeometric differential equations, integer programming, toric geometry, graph theory, e.t.c. (see, e.g. \cite{ABR,ATY, DST, ES,HO,MS, PET1, PET2, ST,VI,VI2}).

A binomial $\textbf{x}^{\textbf{u}}-
\textbf{x}^{\textbf{v}}$ in $I_A$ is called {\it primitive} if
there is no other binomial
 $\textbf{x}^{\textbf{w}}- \textbf{x}^{\textbf{z}}$ in $I_A$,
such that $\textbf{x}^{\textbf{w}}$ divides $
\textbf{x}^{\textbf{u}}$ and $\textbf{x}^{\textbf{z}}$ divides $
\textbf{x}^{\textbf{v}}$. The set of primitive binomials, which is finite, is the
{\it Graver basis} of $I_A$ and is denoted by $Gr_A$.   
 The {\it universal Gr\"{o}bner basis}
 of an ideal $I_A$, is denoted by $\mathcal{U}_A$ and is defined as the union of all reduced Gr\"obner bases $G_\prec$ of $I_A$, as $\prec$ runs over all term orders. Since $I_A$ is generated by binomials, then every reduced Gr\"obner basis of $I_A$ consists of binomials (see, e.g. \cite{ES}). Thus, the universal Gr\"obner basis of $I_A$ is a finite subset of binomials in $I_A$ and it is a Gr\"obner basis for the ideal with respect to all term orders.  By \cite[Proposition 4.11]{ST}, we have that $\mathcal{ U}_A \subseteq  Gr_A$.
  A {\it Markov basis} $M_A$ is a minimal binomial generating set of the toric ideal $I_A$ (its name Markov basis comes from its relation with Markov chains, see \cite[Theorem 3.1]{DST}). The {\it universal Markov basis} of the ideal is denoted by $\mathcal{M}_A$ and is defined as the union of all the Markov bases of the ideal. The elements of $\mathcal M_A$ are called {\it minimal binomials}. Since $A \subseteq \N^n$, then $\mathbb N A$ is {\it pointed} (that is, $\mathbb NA \cap (- \mathbb N A) = \{0\}$). As a consequence  $\mathcal{M}_A \subseteq Gr_A$ and, hence, $\mathcal M_A$ is also a finite set (see \cite[Theorem 2.3]{THO1} and \cite{THO3}).  The Graver basis, the universal Gr\"obner basis and the universal Markov basis are usually called {\it toric bases}.

An ideal $I$ is called {\it robust} if its universal Gr\"obner
basis is a minimal set of generators of the ideal. Even in the context of toric ideals, robustness is a property that has not been fully described. It is known that Lawrence ideals are robust and robustness has also been studied in \cite{BOO2} for toric ideals of graphs, and in \cite{BOO1} for toric ideals which are generated by quadratic binomials.  Some of the interests in studying robustness stems from the fact that they are ideals which are minimally generated by a Gr\"obner basis, see \cite{BOO3, CONCA}. Whenever $I$ is an ideal with a quadratic Gr\"obner basis, then $\mathbb K[x_1,\ldots,x_m]$ is a Koszul algebra. Hence, another interesting feature of robust ideals generated by quadrics, is that they provide examples of Koszul algebras.  Nevertheless, robustness is a rare property and there are very few nontrivial examples of robust ideals. This makes it natural to consider a wider family of ideals that shares many of its good properties. In this paper we study the property of a toric ideal being {\it generalized robust}. An ideal is called generalized robust if its universal Gr\"obner basis is equal to its universal Markov basis. The notion of generalized robustness of a toric ideal was introduced in \cite{CHT} and, as its name indicates, it is a family containing robust toric ideals (see \cite[Corollary 3.5]{CHT}). Since determining or computing the universal Gr\"obner basis of $I$ is a very difficult and computationally demanding problem, it is still difficult to determine whether a toric ideal is generalized robust. The main goal of this paper is to provide several families of generalized robust toric ideals.

The present paper is divided into the following sections.  In Section 2, we collect some basic facts related to the toric bases, that can be found in the literature or are easy consequences of known results. In particular, we give a description of the universal Markov basis (Proposition \ref{pr:minimal}) and study how the toric bases behave with respect to the elimination of variables (Propositions \ref{pr:elimination} and \ref{genrobustfewvariables}). In Theorem \ref{genrobustquarics} we prove that a homogeneous toric ideal is generalized robust and generated by quadrics if and only if its universal Gr\"obner basis only consists of quadrics.

The main results of this paper are divided in two parts.
The first one is presented in Section 3, in which we completely characterize the graphs giving rise to robust and generalized robust toric ideals that both are generated by quadrics. More precisely, we provide the following structural theorem, which summarizes Theorem \ref{quadratic-bipartite} and Theorem \ref{maintheorem} (see also Definition \ref{def:theguy} and Definition \ref{def:necklace}):

\begin{thm1} Let $G$ be a finite, connected and simple graph. The toric ideal associated to $G$ is generalized robust and generated by quadrics if and only if $G$ has at most one non bipartite block which is either a $K_4$, or a $K_{2,\ell}\cup \{e\}$ or a double-$K_{2,(r,s)}$ or a necklace-$K_{2,\ell}$ graph. The bipartite blocks of $G$ are either $K_{2,\ell}$ or cut edges.
\end{thm1}

Interestingly, $\mathbb{K}[x_1,\ldots,x_m]/I_G$ is a Koszul ring for all the graphs $G$ described in this theorem. 

As a direct consequence of this, our second main result is Theorem \ref{maintheorem2}, in which we characterize all graphs whose toric ideal is robust and is generated by quadratic binomials. It should be noticed that this result also completes \cite[Corollary 5.2]{BOO2}, 

\begin{thm1}Let $G$ be a non bipartite graph. The toric ideal associated to $G$ is robust and is generated by quadrics if and only if all the blocks of $G$ are bipartite except one which is either a $K_{2,\ell}\cup \{e\}$ or a double-$K_{2,(r,s)}$ or a necklace-$K_{2,\ell}$ graph. The bipartite blocks of  $G$ are of type $K_{2,\ell}$ or cut edges.
\end{thm1}

The second part of the main results is presented in Section 4, where we work in the framework of toric ideals associated to submonoids of $\mathbb N$. More precisely, given a submonoid $\mS$ of $(\mathbb N, +)$, then it has a unique minimal set of generators $A = \{a_1,\ldots,a_m\}$, and the toric ideal of $\mS$ is defined as $\IS := I_A$. Taking $d := \gcd(a_1,\ldots,a_m)$ and $A' := \{a_1/d,\ldots,a_m/d\}$, then $I_A = I_{A'}$. Hence, one may assume without loss of generality that $A = \{a_1,\ldots,a_m\}$ consists of relatively prime positive integers and, in this case, $\mS$ is called a numerical semigroup (for a detailed study of numerical semigroups we refer the reader to \cite{AG,RA}). 
Since $\IS$ has height $m-1$, we have that $\IS$ is a complete intersection if and only if one of its Markov basis (and, thus, all its Markov bases) consists of $m-1$ binomials.  Complete intersection numerical semigroups have been widely studied in the literature, see, e.g., \cite{AG2,BGS,BGRV,D,DMS,H}. 

Take $\prec$ a monomial order in $\mathbb K[x_1,\ldots,x_m]$ and denote by  ${\rm in}_{\prec}(\IS)$ the corresponding initial ideal of $\IS$. Since  ${\rm ht}({\rm in}_{\prec}(\IS)) = {\rm ht}(\IS) = m-1$, then ${\rm in}_{\prec}(\IS)$ is a complete intersection if it can be generated by $m-1$ monomials. In other words, ${\rm in}_{\prec}(\IS)$ is a complete intersection if and only if the reduced Gr\"obner basis of $\IS$ with respect to $\prec$ consists of $m-1$ binomials. Since Gr\"obner bases are generating sets of the ideal, whenever  ${\rm in}_{\prec}(\IS)$ is a complete intersection for a monomial order $\prec$, then $\IS$ so is. 

Our main results in this section are summarized in the following diagram:

\[ \begin{array}{ccl}  
\mS = \langle a_1, a_2 \rangle & \Longleftrightarrow & \mS \text{ is robust } \\
\Downarrow \\
\mS \text{ has a unique Betti element } & \Longleftrightarrow & \mS \text{ is generalized robust } \\
\Downarrow \\ \mS \text{ is free } & \Longleftrightarrow &   {\rm in}_{\prec}(\IS)  \text{ is a C.I. for some } \prec 
\\  \Downarrow \\  \IS \text{ is a C.I.} \end{array} \]

Theorem \ref{th:existsInitial} states that $\IS$ has a complete intersection initial ideal if and only if $\mS$ is a free numerical semigroup, a family of semingroups studied in \cite{BG, H, RG}. Since in this case, both $\IS$ and the corresponding initial ideal ${\rm in}_\prec(\IS)$ are complete intersections, it turns out that all the Betti numbers in the whole minimal graded free resolution of $\IS$ and ${\rm in}_\prec(\IS)$ coincide, providing examples of {\it robustness of Betti numbers}. This is an interesting phenomenon which is known to occur for robust toric ideals generated by quadrics \cite{BOO1} and also to the ideal of maximal minors of a generic matrix and its Gr\"obner basis with respect to a certain monomial order \cite{CONCA}. In Theorem \ref{th:genrbustnumsemigroup} we determine all toric ideals of numerical semigroup that are generalized robust. It turns out that this property is characterized by a known subfamily of numerical semigroups studied in \cite{GH, GOR, KO}, namely the semigroups with a unique Betti element. As an easy consequence, we get in Corollary \ref{cor:robust} that there are no nontrivial examples of robust toric ideal of a numerical semigroup, that is, $\IS$ is robust if and only if $\IS$ a principal ideal or, in other words, if $\mS$ is a $2$-generated numerical semigroup. 

 Finally, we state some conclusions and formulate some conjectures and open problems concerning robustness, generalized robustness and related properties in toric ideals. These conjectures are supported by some experimental evidence with the computer softwares CoCoA \cite{ABR} and Singular \cite{SINGULAR}.

\section{Remarks on toric bases}
Let $A \subseteq \N^n$ be a finite set of nonzero vectors. In this section we will discuss some known basic facts about the bases associated to the toric ideal $I_A\subseteq \mathbb{K}[{\bf {x}}]:=\mathbb{K}[x_1,\ldots,x_n]$, namely the universal Markov basis $\mathcal M_A$, the universal Gr\"obner basis $\mathcal U_A$, the Graver basis~$Gr_A$ and the set of the circuits $\mathcal C_A$. A binomial $\textbf{x}^{\textbf{u}}- \textbf{x}^{\textbf{v}} \in I_A$ with $\textbf{u} = (u_1,\ldots,u_m),\, \textbf{v} = (v_1,\ldots,v_m) \in \N^m$ is called a {\it circuit} if it has minimal support (with respect to set containment), if $\gcd(\textbf{x}^{\textbf{u}}, \textbf{x}^{\textbf{v}}) = 1$ and the nonzero entries of $\textbf{u} + \textbf{v}$ are relatively prime. 

First of all, it is worth pointing out that, since $\N A$ is pointed, then the graded version of Nakayama's lemma holds. As a consequence, if we consider a set $\{g_1,\ldots,g_r\}$ of $A$-homogeneous polynomials, one has that $I_A = \langle g_1,\ldots,g_r\rangle$ if and only if the cosets of $g_1,\ldots,g_r$ span the $\mathbb K$-vector space $I_A / \langle x_1,\ldots,x_m\rangle \cdot I_A$. As a consequence, any minimal set of $A$-homogeneous generators of $I_A$ has the same cardinality, which is $\mu(I_A) := {\rm dim}_{\mathbb K}(I_A / \langle x_1,\ldots,x_m\rangle \cdot I_A)$.  Moreover, Nakayama's lemma also guarantees that the $A$-degrees appearing in any minimal set of $A$-homogeneous generators are invariant, these values are usually called {\it Betti degrees} of $I_A$. Since every binomial in $I_A$ is $A$-homogeneous, we have the following result.
 \begin{prop1}\label{pr:minimal} The universal Markov basis $\mathcal M_A$ is the set of binomials in $I_A$ that do not belong to $ \langle x_1,\ldots,x_m\rangle \cdot I_A$.
 \end{prop1}

Moreover, the following inclusions of toric bases hold:
 
 \begin{thm1} \label{pr:inclusions} \cite[Proposition 4.11]{ST}, \cite[Theorem 2.3]{THO1} For any toric ideal it holds
 $$\mathcal C_A \subseteq \mathcal U_A \subseteq Gr_A.$$
 Moreover, since $A$ defines a pointed semigroup, then $\mathcal M_A \subseteq Gr_A.$
 \end{thm1}
 
In the forthcoming we use many times how the toric bases behave with respect to an elimination of variables, these are summarized in the following result (see, e.g., \cite[Proposition 4.13]{ST}).
 
 \begin{prop1}\label{pr:elimination} Let $A' \subseteq A$ and denote $\mathbb K[\x_{A'}] :=  \mathbb K[x_i \, \vert \, a_i \in A']$, then 
\begin{itemize}
\item[(a)] $I_{A'} = I_A \cap \mathbb K[\x_{A'}]$.
\item[(b)] $\mathcal C_{A'} = \mathcal C_A \cap \mathbb K[\x_{A'}]$.
\item[(c)] $\mathcal U_{A'} = \mathcal U_A \cap \mathbb K[\x_{A'}]$.
\item[(d)] $Gr_{A'} = Gr_A \cap \mathbb K[\x_{A'}]$.
\item[(e)] $\mathcal M_{A'} \supseteq \mathcal M_A \cap \mathbb K[\x_{A'}]$.

\end{itemize}
 \end{prop1}

In general we do not have equality in Proposition \ref{pr:elimination}.(e). For example, for $A' = \{a_1,a_2\} \subseteq A = \{a_1,a_2,a_3\} \subseteq \N$ with $a_1 = 4, a_2 = 5$ and $a_3 = 6$, we have that $\mathcal M_A = \{x_1^3 - x_3^2, x_2^2 - x_1 x_3\}$, and hence $\mathcal M_A \cap \mathbb K[x_1,x_2] = \emptyset$, whereas  $\mathcal M_{A'} = \{x_1^5 - x_2^4\}$.

In \cite[Proposition 2.5]{BOO1} the authors proved that robustness is preserved under an elimination of variables, that is, if an ideal $I_A$ is robust and $A' \subseteq A$, then the ideal $I_{A'}$ is also robust. We do not know if generalized robustness is preserved under elimination of variables, see Question \ref{question}. Nevethelesss, the following results hold.

\begin{prop1}\label{genrobustfewvariables} Let $A \subseteq \N^n$ be a finite set of nonzero vectors and $A' \subseteq A$. Then, 
\begin{itemize}
\item[(a)] If $\mathcal U_A \subseteq \mathcal M_A$, then $\mathcal U_{A'} \subseteq \mathcal M_{A'}$.
\item[(b)] If $\mathcal C_A \subseteq \mathcal M_A$, then $\mathcal C_{A'} \subseteq \mathcal M_{A'}$.
\item[(c)] If $Gr_A =  \mathcal M_A$, then $Gr_{A'} = \mathcal M_{A'}$.
\end{itemize} 
\end{prop1}
\begin{proof} 
(a) If $\mathcal{U}_A \subseteq \mathcal{M}_A$, by Proposition \ref{pr:elimination} we have that
$\mathcal U_{A'} =  \mathcal U_A \cap \mathbb K[\x_{A'}] \subseteq \mathcal{M}_A \cap \mathbb K[\x_{A'}] \subseteq \mathcal M_{A'}.$
The proof of (b) is analogue to the one of (a). 

(c) If $Gr_A = \mathcal{M}_A$, by Proposition \ref{pr:elimination} we have that
$Gr_{A'} =  Gr_A \cap \mathbb K[\x_{A'}] = \mathcal{M}_A \cap \mathbb K[\x_{A'}] \subseteq \mathcal M_{A'},$ the inclusion $\mathcal M_{A'} \subseteq Gr_{A'}$ follows from Proposition \ref{pr:inclusions}.
\end{proof}

The next theorem gives a nice property which motivates the study of generalized robust toric ideals which are generated by quadrics.

\begin{thm1}\label{genrobustquarics} A homogeneous toric ideal $I_A \subseteq \mathbb{K}[x_1,\ldots,x_m]$ is generalized robust and generated by quadrics if and only if the universal Gr\"obner basis of $I_A$ only consists of quadrics.
\end{thm1}
\begin{proof}$(\Longrightarrow)$ Since $I_A$ is generated by quadrics, then its universal Markov basis $\mathcal M_A$ consists of quadrics. 
The result follows from the definition of a generalized robust toric ideal. 

$(\Longleftarrow)$ Since the universal Gr\"obner basis is a set of generators, then $I_A$ is generated by quadrics. Let us see now that $\mathcal M_A = \mathcal U_A$. 

Let $f \in \mathcal U_A$, then $f \in I_A$ is a quadric and $f \notin \langle x_1,\ldots,x_m\rangle \cdot I_A$ (since the elements of $\langle x_1,\ldots,x_m\rangle \cdot I_A$ have degree at least three). By the graded version of Nakayama's lemma it follows that $f \in \mathcal M_A$ and thus $\mathcal U_A\subseteq \mathcal M_A$. Conversely, take $f \in \mathcal M_A$. Since the ideal $I_A$ is generated by quadrics and due to the fact that there exists a Markov basis of the ideal which consists of quadrics, by Nakayama's lemma it follows that every Markov basis of $I_A$ consists of quadrics and thus $f$ is a quadratic binomial. After reindexing the variables and considering $-f$ if necessary, we have that either $f = x_3 x_4 - x_1 x_2$, or $f = x_3^2 - x_1x_2$, otherwise the ideal is not prime. In both cases the binomial $f$ is in the reduced Gr\"obner basis of $I_A$ with respect to the lexicographic order with $x_m > \cdots > x_1$. Thus, $f \in \mathcal U_A$ and therefore $\mathcal M_A\subseteq \mathcal U_A$.
\end{proof}

\section{Quadratic robust and generalized robust toric ideals of graphs}\label{sec:graph}

\subsection{Preliminaries}
 In this section we study the robustness and generalized robustness property for toric ideals of graphs which are generated by quadratic binomials. In the rest of the present section, we consider finite, simple and connected graphs. Let $G$ be a graph with vertices
$V(G)=\{v_{1},\ldots,v_{n}\}$ and edges $E(G)=\{e_{1},\ldots,e_{m}\}$.
Let $\mathbb{K}[e_{1},\ldots,e_{m}]$
be the polynomial ring in the $m$ variables $e_{1},\ldots,e_{m}$ over a field $\mathbb{K}$.  We
 associate each edge $e=(v_{i},v_{j})\in E(G)$ with the element
$a_{e}=v_{i}+v_{j}$ in the free abelian group $ \mathbb{Z}^n $
with basis the set of the vertices
of $G$, i.e. each vertex $v_j\in V(G)$ is associated with the vector $(0,\ldots,0,1,0,\ldots,0)\in \mathbb{N}^n$, where the 
nonzero component is in the $j$ position. We denote by $I_G$ the toric ideal $I_{A_{G}}$ in
$\mathbb{K}[e_{1},\ldots,e_{m}]$, where  $A_{G}=\{a_{e}\ | \ e\in E(G)\}\subseteq \mathbb{N}^n $.
Toric ideals of graphs are homogeneous prime ideals and many of their algebraic properties can be described in terms
of the underlying graph. For example, $I_G$ has height \begin{equation}\label{eq:heightgraphideal}
{\rm ht}(I_G)= |E(G)| -  |V(G)| + b(G),
\end{equation} where  $b(G)$ denotes the number of connected components of the graph $G$ that are bipartite (see, e.g., \cite{V}). In particular, if $G$ is connected, then $b(G) = 1$ if $G$ is bipartite and $b(G) = 0$ otherwise (see, e.g., \cite{GRV}). 
Also, binomial generating sets of $I_G$ can be described in terms of some walks in the graph. To present this result, we first recall some basic elements from graph theory (for unexplained terminology and basics on graphs we refer to \cite{BM}).

A {\it walk} connecting $u \in V(G)$ and
$u' \in V(G)$ is a finite sequence of vertices of the graph 
$w=(u = u_0, u_1, \ldots, u_{\ell-1}, u_\ell = u')$,
with each $e_{i_j}=(u_{j-1},u_{j})\in E(G)$, for $j=1,\ldots,\ell$. 
The {\it length}
of the walk $w$ is the number $\ell$ of its edges. An
{\it even} (respectively {\it odd}) {\it walk} is a walk of even (respectively odd) length.
A walk
$w=(u_0,u_1\ldots,u_{\ell-1},u_{\ell})$
is called {\it closed} if $u_{0}=u_{\ell}$. A {\it cycle}
is a closed walk
$(u_0,u_{1},\ldots,u_{\ell-1},u_{\ell})$ with
$u_{k}\neq u_{j},$ for every $ 1\leq k < j \leq \ell$. 

Consider an even closed walk $w=(u_{0},u_{1},u_{2},\ldots,u_{2s-1},u_{2s} = u_0)$ of length $2s$ with $e_{i_j}=(u_{j-1},u_{j})\in E(G)$, for $j=1,\ldots,2s$. The binomial \begin{equation}\label{eq:bingen}
B_w = e_{i_1}e_{i_3} \cdots e_{i_{2s-1}}  - e_{i_2} e_{i_4} \cdots e_{i_{2s}}
\end{equation} 
belongs to the toric ideal $I_G$. Actually, Villarreal proved  in \cite{VI2} that 
\begin{equation}\label{eq:idealgraph}
I_G = \langle B_w  \ \vert \ w \text{ is an even closed walk}\rangle,
\end{equation} 
that is, the toric ideal $I_G$ is generated by the binomials corresponing to even closed walks of the graph $G$.

In \cite{OH}, Ohsugi and Hibi gave the following combinatorial criterion for the toric ideal of a graph $G$ to be generated by quadrics:

\begin{thm1}\cite[Theorem 1.2]{OH}\label{idealbyquadratic} Let $G$ be a finite connected simple graph. Then, the toric ideal $I_G$ of $G$ is generated by quadrics if and only if the following conditions are satisfied:
\begin{itemize}
\item[i)] if $c$ is an even cycle of $G$ of length $\geq 6$, then either $c$ has an even chord or $c$ has three odd chords $e, e', e''$ such that $e$ and $e'$ cross in $c$,
\item[ii)] if $c_1$ and $c_2$ are odd chordless cycles of $G$ having exactly one common vertex, then there exists a bridge between them, 
\item[iii)] if $c_1$ and $c_2$ are odd chordless cycles of $G$ having no common vertex, then there exist at least two bridges between $c_1$ and $c_2$.
\end{itemize}
\end{thm1}

Also, the only even closed walks of length four in a simple graph are cycles. Hence, whenever $B_w$ is a quadric, then $w$ is a cycle of length four. Thus, we have the following.

\begin{cor1}\label{quadratic=allminimal=4cycles} Let $G$ be a graph. If the toric ideal $I_G$ is generated by quadrics, then all its minimal binomials are of the form $B_w$, where $w$ is a cycle of length four.
\end{cor1}

The Graver basis and the universal Markov basis of the toric ideal of a graph $G$, which we denote by $Gr_G$ and $\mathcal M_G$ correspondingly, were described in \cite{TT1}, Theorem 3.2 and Theorem 4.13 correspondingly, while the universal Gr\"obner basis, which we denote by $\mathcal U_G$, was described in \cite[Theorem 3.4]{TT2}. For the sake of brevity we refer the reader to the above articles. A necessary and sufficient characterization for generalized robust toric ideals of graphs was given in \cite{CHT}:

\begin{thm1}\cite[Theorem 3.4]{CHT}\label{M=Gr} Let $G$ be a graph and let $I_G$ be its corresponding toric ideal. The ideal $I_G$ is
generalized robust if and only if $\mathcal{M}_G=Gr_G$.
\end{thm1}

For general toric ideals, it may happen that the universal Markov basis is not contained into the universal Gr\"obner basis (see Section \ref{sec:conclusion} for an example). Nevertheless, in the context of toric ideals of graphs we have that $\mathcal M_G \subseteq \mathcal U_G$ for any graph $G$ (see \cite[Proposition 3.3]{CHT}).  This fact together with Proposition \ref{genrobustfewvariables} yields that the generalized robustness property of a graph $G$ is a hereditary property, in the sense that it holds also for any subgraph of $G$.

\begin{cor1}\label{hereditary} Let $H$ be a subgraph of a graph $G$. If the ideal $I_G$ is generalized robust, then $I_{H}$ is generalized robust.
\end{cor1}

\subsection{Quadratic generalized robust graphs$;$ the bipartite case}

We state some properties of a generalized robust toric ideal of a graph, which stem from directly of the results in \cite{TT1, CHT, TT2} and will be useful for us in the sequel. By {\it chordless graph} we mean a graph in which every cycle has no chords.

\begin{cor1}\label{genrobustbipartite-chordless+4} Let $G$ be a bipartite graph.
\begin{itemize}
\item[$\alpha$)] The ideal $I_G$ is generalized robust if and only if the graph $G$ is chordless. 
\item[$\beta$)] The ideal $I_G$ is generalized robust and generated by quadrics if and only if all the cycles of the graph $G$ have length four.
\end{itemize} 
\end{cor1}
\begin{proof}
$(\alpha)$ In \cite[Proposition 4.3]{CHT} the author proved that an even cycle in a  graph of a generalized robust toric ideal has only odd chords (if it has). Since the graph $G$ is bipartite it follows that it is chordless. Conversely, it is known that $\mathcal{M}_G\subseteq Gr_G$ and let $B_w\in Gr_G$. The graph is bipartite, thus the walk $w$ is an even cycle, see \cite[Theorem 3.2]{TT1}. Since the graph is chordless, the cycle $w$ is chordless and therefore $B_w\in \mathcal{M}_G$, see \cite[Theorem 4.13]{TT1}. The result follows from Theorem \ref{M=Gr}.

$(\beta)$ It follows from the previous argument ($\alpha$) and from Theorem \ref{idealbyquadratic}.

\end{proof}

\begin{figure}[h]
\begin{center}
\includegraphics[scale=0.7]{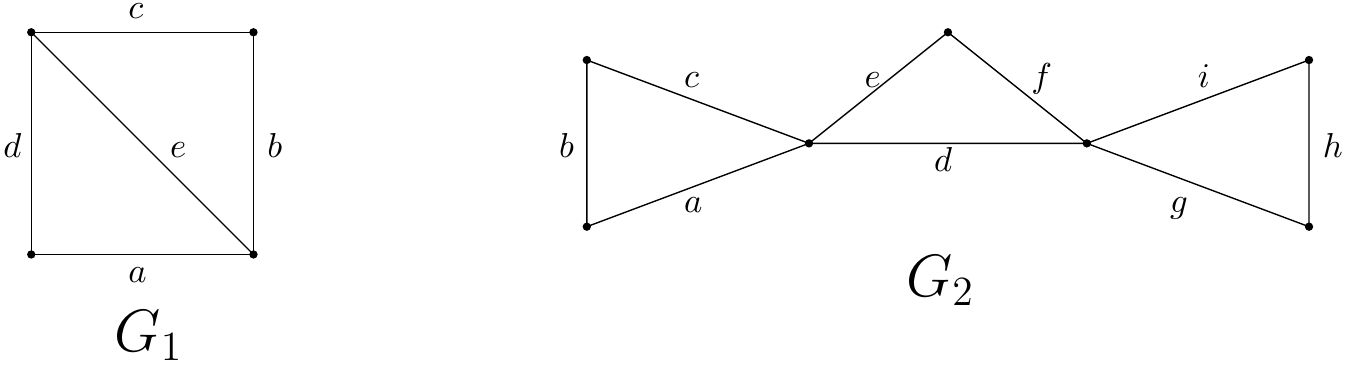}
\caption{Corollary \ref{genrobustbipartite-chordless+4}.($\alpha$) does not hold for non bipartite graphs.}
\label{figure-chordless}
\end{center}
\end{figure}

We remark that in the case of non bipartite graphs, none of the two of the implications of Corollary \ref{genrobustbipartite-chordless+4} $(\alpha)$ are true. For example in Figure \ref{figure-chordless}, we present a non bipartite graph $G_1$ which is not chordless and whose toric ideal is generalized robust since it is principal with $I_{G_1}=\langle ac-bd\rangle$, see \cite[Theorem 4.13]{TT1}. Also, in Figure \ref{figure-chordless}, we present a chordless non bipartite graph $G_2$ whose corresponding toric ideal is not generalized robust. Indeed, one can verify that it not satisfies the condition in Theorem \ref{M=Gr}. More precisely, the binomial $B_w=acf^2h-be^2ig$ is in the Graver basis of the ideal and is not minimal since the corresponding walk $w=(a,b,c,e,f,i,h,g,f,e)$ has a bridge $d$, see \cite[Theorems 3.2 and 4.13]{TT1}.

By Corollary \ref{genrobustbipartite-chordless+4}, the bipartite graphs yielding generalized robust toric ideals of graphs which are generated by quadrics are exactly the graphs whose cycles have all length $4$. In Theorem \ref{quadratic-bipartite} we describe precisely these graphs. Before proceeding with its statement and proof, we introduce some definitions and notation.

\begin{figure}[h]
\begin{center}
\includegraphics{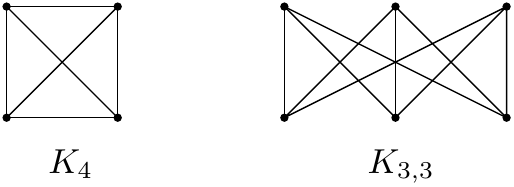}
\caption{The complete graphs $K_4$ and $K_{3,3}$}
\label{k4}
\end{center}
\end{figure}

We denote by $K_n$ the complete graph on $n$ vertices and by $K_{r,s}$ the complete bipartite graph with partitions of sizes $r \in \mathbb Z^+$ and $s \in \mathbb Z^+$. In Figure \ref{k4} we see the complete graph on four vertices $K_4$ and the complete bipartite graph $K_{3,3}$.

A cut edge (respectively cut vertex) is an edge (respectively vertex) of
the graph whose removal increases the number of connected
components of the remaining subgraph.  A graph is called biconnected if it is connected and does not contain a cut
vertex. A block is a maximal biconnected subgraph of a 
graph. For a given graph $G$ and a set $S\subseteq V(G)$, we denote by $[S]$ the corresponding induced subgraph of $G$, that is,  the graph with vertex set $S$ and whose edge set consists of the edges in $E(G)$ having both endpoints in $S$. 

We have the following combinatorial lemma.

\begin{lem1}\label{lm:4cycles}The only biconnected graphs whose cycles have all length $4$ are cut edges and complete bipartite graphs $K_{2,\ell}$ with $\ell \geq 2$. 
\end{lem1}
\begin{proof} If the graph has no cycles, then it is biconnected if and only if it is a cut edge. The graph $K_{2,\ell}$  has only cycles of length $4$ for all $\ell \geq 2$. Consider now a non-acyclic biconnected graph $B$ whose cycles all have length $4$. We have that $B$ is bipartite and we denote by $U$ and $V$ the bipartition of $V(B)$. Take $u \in U$ and $v \in V$, since $B$ is biconnected, there are two disjoint paths from $u$ to $v$. Moreover, $B$ has only cycles of length four, so one of these paths has length $1$ and, hence, $\{u,v\}$ is an edge of $B$. As a consequence $B$ is a complete bipartite graph $B = K_{t,\ell}$ for some $2 \leq t \leq  \ell$. If $t \geq 3$, then $B$ has a $6$-cycle, a contradiction. Thus, we conclude that $B = K_{2,\ell}$ for some $\ell \geq 2$. 
\end{proof}

As a consequence of this lemma and Corollary \ref{genrobustbipartite-chordless+4} we get the following. 

\begin{thm1}\label{quadratic-bipartite} Let $G$ be a bipartite graph. The ideal $I_G$ is generalized robust and generated by quadrics if and only if all the blocks of $G$ are $K_{2,\ell}$ or cut edges, for some $\ell \geq 2$.
\end{thm1}
\begin{proof}
($\Longrightarrow$) Let $I_G$ be a generalized robust toric ideal which is generated by quadratic binomials and let $B$ be one of the blocks of $G$. By Corollary \ref{hereditary}  and Lemma \ref{lm:4cycles} we have that $B$ is a cut edge or a $K_{2,\ell}$, for some $\ell \geq 2$.

($\Longleftarrow$) The graph $G$ is consists of blocks which are either $K_{2,\ell}$ or cut edges. Therefore all the cycles of $G$ have length four. By Corollary \ref{genrobustbipartite-chordless+4} ($\beta$) it follows that the ideal is generalized robust and is generated by quadrics.
\end{proof}

\begin{figure}[h]
\begin{center}
\includegraphics[scale=0.7]{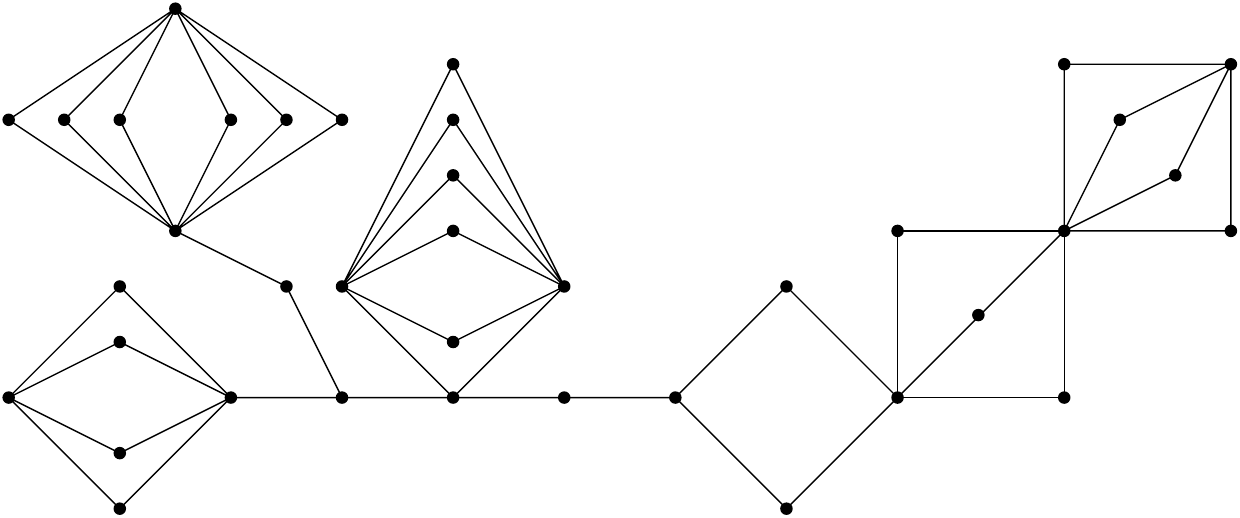}
\caption{A bipartite graph $G$ such that the ideal $I_G$ is a quadratic generalized robust}
\label{biprobust}
\end{center}
\end{figure}

In Figure \ref{biprobust} we present an example of a bipartite graph $G$ whose corresponding toric ideal is generalized robust and is generated by quadrics.

\subsection{Quadratic generalized robust graphs$;$ the general case}

We are moving on to the general case of non bipartite graphs. As we can see in the next lemma, when the toric ideal of a graph is generated by quadratic binomials, then it has at most one non bipartite block.

\begin{lem1} \label{quadraticatmostonenonbipar} Let $G$ be a simple connected graph such that the corresponding toric ideal $I_G$ is generated by quadrics. Then the graph $G$ has at most one non bipartite block.
\end{lem1}
\begin{proof} Let $G$ be a connected graph with at least two non bipartite blocks and let them be $B_1$ and $B_2$.  Let $c_1$ and $c_2$ be two chordless odd cycles of the blocks $B_1$ and $B_2$ correspondingly. We separate the proof in two cases: either the cycles have (exactly) one common vertex or they are vertex disjoint. 

In the first case, by Theorem \ref{idealbyquadratic}.(ii) there exists a bridge between the cycles $c_1,c_2$, but this contradicts the fact that the cycles belong to different blocks. In the second case, Theorem \ref{idealbyquadratic}.(iii) guarantees that there exist $e_1=(x_1,y_1)$ and $e_2=(x_2,y_2)$ two bridges between $c_1,c_2$, where $x_1,x_2\in V(c_1)$ and $y_1,y_2\in V(c_2)$. There are two cases, either $x_1\neq x_2$ and $y_1\neq y_2$ or $x_1=x_2$ (similarly if $y_1=y_2$). The first case is not possible because $c_1,c_2$ are in different blocks. The second case yields two odd cycles in different blocks with one common vertex (the vertex $x_1=x_2$), and we already discussed that this is not possible. 
\end{proof}

In order to give the structural characterization of a graph $G$ such that $I_G$  is generalized robust and generated by quadrics, we need to introduce the notions of the {\it double-$K_{2,(r,s)}$ graph} and the {\it necklace-$K_{2,\ell}$ graph}. We remind that a {\it 2-clique sum} of the graphs $G_1$ and $G_2$ is obtained by identifying an edge $e_1$ of $G_1$ and an edge $e_2$ of $G_2$.

\begin{def1}\label{def:theguy}We consider the non bipartite graph $G_1=K_{2,r}\cup \{e\}$, where $e$ is an edge connecting two vertices of $K_{2,r}$ and the bipartite graph $G_2=K_{2,s}$, where $r,s\geq2$. A graph $G$ is called a double-$K_{2,(r,s)}$ if it is a 2-clique sum of the graphs $G_1$ and $G_2$ obtained by identifying the edge $e$ of $G_1$ with any edge of $G_2$.
\end{def1}

\begin{figure}[h]
\begin{center}
\includegraphics[scale=0.6]{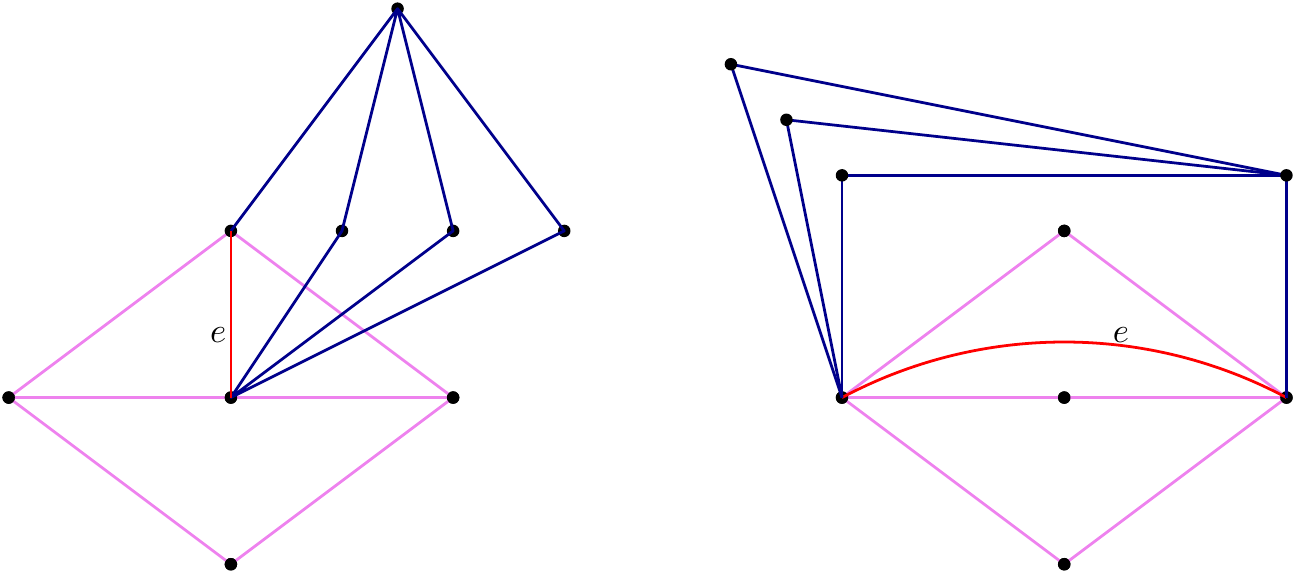}
\caption{The two non-isomorphic double-$K_{2,(3,4)}$ graphs}
\label{theguys}
\end{center}
\end{figure}

For example, in Figure \ref{theguys} we present the two non-isomorphic double-$K_{2,(3,4)}$ graphs.

\begin{prop1}  \label{pr:theguy} If $G$ is a double-$K_{2,(r,s)}$, then $I_G$ is a generalized robust toric ideal generated by quadrics.
\end{prop1}
\begin{proof} Let $G$ be a double-$K_{2,(r,s)}$ graph and we consider the graph $G'$ which consists of two connected components, the $K_{2,r}$ and the $K_{2,s}$. We observe that $I_{G'} \subseteq I_G$ since every even closed walk in $G'$ is in $G$. Since both connected components of $G'$ are bipartite, we have that $b(G') = 2$ and, by~(\ref{eq:heightgraphideal}):
 \begin{eqnarray}
{\rm ht}(I_{G'}) &=&  |E(G')| - |V(G')| + 2 = 2(r+s) - (r+s+4) +2 \nonumber \\
&=& 2(r+s) - (r+s+2) = r+s-2 \nonumber\\
&=&  |E(G)| - |V(G)| = {\rm ht}(I_G). \nonumber
\end{eqnarray}

So both $I_G$ and $I_{G'}$ are prime ideals of the same height and $I_{G'} \subseteq I_G$, hence $I_G = I_{G'}$ and the proof follows from Theorem \ref{quadratic-bipartite}.  
\end{proof}

In \cite{TT3} the authors define the necklace graph as a graph which comes from identifying two vertices at odd distance of a chain of bipartite blocks. Following the same structure, we define the necklace-$K_{2,\ell}$ graphs. Let $T_G$ be the {\it block tree} of a graph $G$, that is, the bipartite graph with bi-partition $( \mathbb{B},S)$ where $\mathbb{B}$ is the set of the blocks of $G$ and $S$ is the set of the cut vertices of $G$, such that $(B,u)$ is an edge if and only if $u\in B$.  A chain of bipartite blocks is a graph $G$ such that its block tree $T_G$ is a path.

\begin{def1}\label{def:necklace} Let $R$ be a bipartite graph consisting of a chain of (bipartite) blocks $B_1,\ldots,B_k$ where $k \geq 2$ and either $B_i=K_{2,n_i}$ for some $n_i \geq 2$, or $B_i$ is a cut edge of $R$, for $i=1,\ldots,k$. Let $x_1 \in V(B_1)$ and $x_2 \in V(B_k)$ be two non-adjacent vertices of $R$ at odd distance which are not cut vertices of $R$. We define a necklace-$K_{2,\ell}$ graph as the graph $G$ obtained after identifying the vertices $x_1$ and $x_2$. That is, the graph on the vertex set $$V(G) = (V(R) \setminus \{x_1,x_2\}) \cup \{x\}$$ and edges $$E(G) = E(R \setminus \{x_1,x_2\}) \cup \{\{u,x\} \, \vert \, \textrm{ either } \{u,x_1\} \in E(R) \text{ or } \{u,x_2\} \in E(R)\}.$$
\end{def1}

\begin{figure}[h]\label{fig:necklace}
\includegraphics[scale=.9]{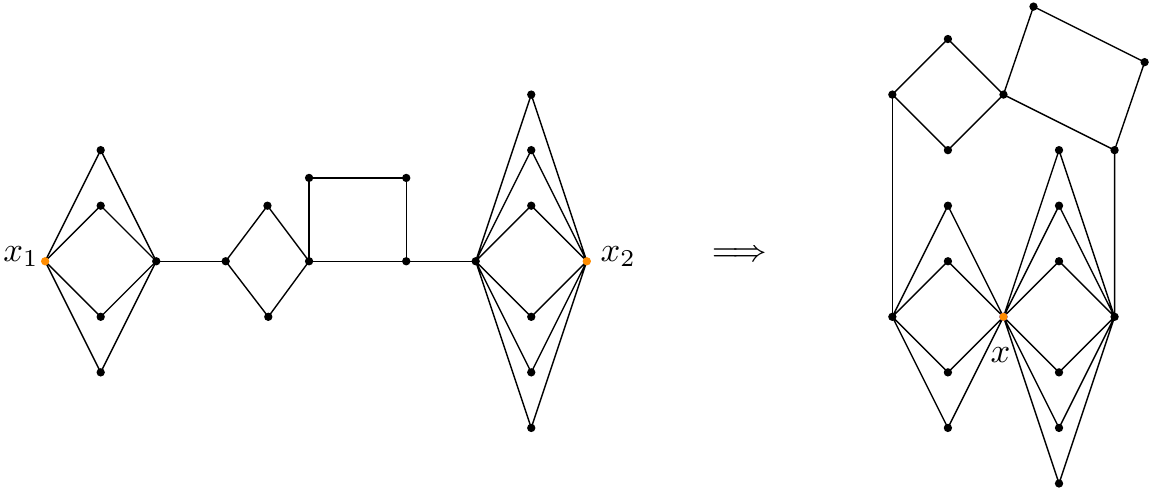}
\caption{The construction of a necklace-$K_{2,\ell}$ graph}
\end{figure}

\begin{prop1}\label{pr:necklace} If $G$ is a necklace-$K_{2,\ell}$, then $I_G$ is a generalized robust toric ideal generated by quadrics.
\end{prop1}
\begin{proof}
Take $R$ as in the definition of necklace-$K_{2,\ell}$. Since every walk in $R$ is a walk in $G$ we have that $I_R \subseteq I_G$. Moreover, the graph $R$ is bipartite and therefore ${\rm ht}(I_R) =|E(R)| - |V(R)| +1$, while the graph $G$ is not bipartite, that is ${\rm ht}(I_G) = |E(G)| - |V(G)|$. By (\ref{eq:heightgraphideal}) it follows that 
\begin{eqnarray}
{\rm ht}(I_{G}) &=&  |E(G)| - |V(G)| =  |E(R)| - \left( |V(R)|- 1\right) \nonumber \\
&=& |E(R)| - |V(R)| +1 = {\rm ht}(I_R). \nonumber
\end{eqnarray}
Hence $I_R = I_G$ and the proof follows from Theorem \ref{quadratic-bipartite}.  
\end{proof}

Next, we state the main result of this section in which we give a complete characterization of the graphs $G$ such that $I_G$ is a generalized robust toric ideal generated by quadratic binomials. By Theorem \ref{genrobustquarics}, this class coincides with the graphs $G$ such that $\mathcal U_G$ consists of quadrics.

\begin{thm1} \label{maintheorem}Let $G$ be a non bipartite graph. The ideal is generalized robust and is generated by quadrics if and only if all the blocks of $G$ are bipartite except one, which is either a $K_4$ or a $K_{2,\ell}\cup \{e\}$ or a double-$K_{2,(r,s)}$ or a necklace-$K_{2,\ell}$ graph. Every bipartite block of the graph $G$ are of type $K_{2,\ell}$ for some $\ell$ or cut edges.
\end{thm1}

In the poof of this result we repeatedly use the following lemma, which is a consequence of \cite[Corollary 3.3]{TT1}.

\begin{lem1}\label{lm:shortcycles} Let $G$ be a connected graph. If $I_G$ is generalized robust and generated by quadrics, then:
\begin{itemize} \item[(a)] it has no even cycles of length $\geq 6$,
\item[(b)] it has no two edge disjoint odd cycles.
\end{itemize}
\end{lem1} 
\begin{proof}(a) Let $w$ be an even cycle of $G$. By \cite[Corollary 3.3]{TT1} we have that $B_w \in Gr_G$. Since $I_G$ is generalized robust, by Theorem \ref{M=Gr}, we have that $B_w \in \mathcal M_G$. But since $I_G$ is generated by quadrics, then $B_w$ is quadric and we conclude that $w$ is a $4$-cycle. 

(b) Assume that there are two edge disjoint odd cycles. By \cite[Corollary 3.3]{TT1}, there exists an even closed walk $w$ (which consists of the above two odd cycles and a walk connecting them) such that $B_w \in Gr_G$. Since the ideal is generalized robust, by Theorem \ref{M=Gr} we have that $B_w$ is also minimal and has degree $\geq 3$, a contradiction due to the fact that the ideal is generated by quadrics.
\end{proof}

\noindent {\it Proof of Theorem \ref{maintheorem}.}  ($\Longrightarrow$) Assume that $I_G$ is generalized robust and generated by quadrics. By Lemma \ref{quadraticatmostonenonbipar}, the graph $G$ has exactly one non bipartite block which we denote by~$B$. By Corollary \ref{hereditary} and Theorem \ref{quadratic-bipartite} the rest of the blocks are of the form $K_{2,\ell}$ or cut edges. We separate two cases: either $(i)$ there exists an edge $e\in E(B)$ such that the graph $B\setminus \{e\}$ is bipartite, or $(ii)$ for every edge $e\in E(B)$ the graph $B\setminus \{e\}$ is not bipartite; where $B\setminus \{e\}$  denotes the graph with vertex set $V(B)$ and edge set $E(B)-\{e\}$.
 
 $(i)$ Let $e \in E(B)$ such that $B\setminus \{e\}$ is bipartite. Combining  Theorem \ref{quadratic-bipartite} and  Corollary \ref{hereditary}, it follows that the blocks of $B\setminus \{e\}$ are $K_{2,\ell}$ or cut edges. If $B \setminus \{e\}$ is still biconnected, then $B$ is a $K_{2,\ell}\cup \{e\}$. If $B \setminus \{e\}$ is not biconnected, then it is a chain of blocks, where each block is either a $K_{2,\ell}$ or a cut-edge. We get thus that $B$ is a necklace-$K_{2,\ell}$ graph. 

  $(ii)$ Assume that for every edge $e\in E(B)$ the graph $B\setminus \{e\}$ is not bipartite.
  
  { \it Claim 1:} Every edge of $B$ belongs to a cycle of length four of $B$.
  
  Suppose not and let $\epsilon$ be an edge of $B$ that does not belong to a cycle of length four of $B$. Since $I_B$ is generated by quadrics and every $4$-cycle of $B$ is in $B - \{\epsilon\}$, by Corollary \ref{quadratic=allminimal=4cycles},  it follows that $I_B\subseteq I_{B\setminus \{\epsilon\}}$.  Obviously we have that $I_{B\setminus \{\epsilon\}}\subseteq I_B$ and therefore $I_{B\setminus \{\epsilon\}}= I_B$. But the graph $B\setminus \{\epsilon\}$ is not bipartite and connected, then 
 \begin{eqnarray}
{\rm ht}(I_{B\setminus \{\epsilon\}}) &=& |E(B\setminus \{\epsilon\})|-|V(B\setminus \{\epsilon\})|  \nonumber \\
&=& |E(B)|-1-|V(B)|=ht(I_B)-1, \nonumber
\end{eqnarray}  a contradiction.
  
  Let $G_1,\ldots,G_k$ be the maximal subgraphs (with respect to the inclusion) of type $K_{2,\ell}$ of $B$. By (Claim 1) we have that every edge of the block $B$ belongs to a $K_{2,\ell}$, thus we have that $E(G_1)\cup \ldots \cup E(G_k)=E(B)$.
  
{ \it Claim 2:} $E(G_i)\cap E(G_j)=\emptyset$ for all $1 \leq i <  j \leq k.$
 
Otherwise, suppose that there exist distinct $i,j\in\{1,\ldots,k\}$ such that  $G_i$ and $G_j$ have at least one edge in common. Let $G_i = K_{2,r}$ and $G_j = K_{2,s}$ with $r,s \geq 2$, and denote by $\{v_1,v_2\}, \{w_1,w_2,\ldots,w_r\}$ and $\{u_1,u_2\}, \{x_1,x_2,\ldots,x_s\}$ their corresponding bipartitions. Let $e$ be the common edge of $G_i, G_j$ and without loss of generality we suppose that $e=(v_1,w_1)=(u_1,x_1)$ (otherwise we rename the vertices). We note that the vertex $u_2 \notin V(G_i)$. Indeed,  $u_2=v_2$ contradicts the maximality of $G_j$, and $u_2=w_2$ implies that the odd cycle $(v_1,w_1,u_2)$ belongs to the bipartite graph $G_i$, a contradiction. Similarly, we note that the vertex $x_2\notin V(G_i)$.  We conclude that we have a length $6$ cycle $(v_1=u_1,x_2,u_2,w_1=x_1,v_2,w_2,v_1)$ in $B$, a contradiction  to Lemma \ref{lm:shortcycles}, and (Claim 2) is proved.

We denote by $[G_1],\ldots,[G_k]$ the induced subgraphs with vertices $V(G_1),\ldots,V(G_k)$ correspondingly. We split the proof in two cases: either $(ii_a)$ the graphs $[G_1],\ldots,[G_k]$ are bipartite or $(ii_b)$ there exists $i\in\{1,\ldots,k\}$ such that the graph $[G_i]$ is non bipartite.

$(ii_a)$ First of all, we remark that $|V(G_i)\cap V(G_j)| \ \leq 1$ for all $1 \leq i < j \leq k$. 

Otherwise, suppose that there exist distinct $i,j\in\{1,\ldots,k\}$ such that  $G_i$ and $G_j$ have at least two vertices in common. 
Let $G_i = K_{2,r}$ and $G_j = K_{2,s}$ with $r,s \geq 2$, and denote by $\{v_1,v_2\}, \{w_1,w_2,\ldots,w_r\}$ and $\{u_1,u_2\}, \{x_1,x_2,\ldots,x_s\}$ their corresponding bipartitions. Since both $G_i$ and $G_j$ are complete bipartite graphs and they do not share edges by (Claim 2), then the two common vertices are not adjacent. 
By the maximality of $G_i$ and $G_j$, the common vertices have to be $x_i=w_{i'}$ and $x_j=w_{j'}$, for some $i,j\in\{1,\ldots,r\}$ and $i',j'\in\{1,\ldots,s\}$ and $r,s > 2$. Then a cycle of length 6 arises$;$ namely the cycle $(x_i=w_{i'}, u_1,x_j=w_{j'},v_1,x_k,v_2,x_i=w_{i'})$ with $k \in \{1,\ldots,s\} \setminus \{i,j\}$, a contradiction to Lemma \ref{lm:4cycles}.  

Consider now $G_1,G_2$ two maximal subgraphs of type $K_{2,\ell}$ with one vertex in common (there are such subgraphs since $B$ is biconnected) and let $\{x\}$ be the common vertex of $G_1,G_2$ (the vertex which we discussed above). Take $G'$ the graph with vertex set $$V(G')=V(B\setminus \{x\}) \cup \{x_1,x_2\}$$ and edge set $$E(G')=E(B\setminus \{x\})\cup \{\{x_1,u\} : \{u,x\}\in E(G_1)\}\cup \{\{x_2,v\} : \{v,x\}\in E(B\setminus G_1)\}.$$ 

By (Claim 1) and (Claim 2) we know that $E(B)=E(G_1)\sqcup \ldots \sqcup E(G_k)$, then $|E(G')| = |E(B)|$.
 As a consequence, we have that 
\begin{equation} \label{eq:ht1} {\rm ht}(I_{G'})= \begin{cases} |E(B)|- (|V(B)|+1) +1,\ \textrm{if}\ G' \ \textrm{is bipartite} \\
|E(B)| - (|V(B)| +1) ,\ \textrm{if}\ G' \ \textrm{is not bipartite} \end{cases} \end{equation}

while \begin{equation} \label{eq:ht2} {\rm ht}(I_B)=|E(B)| - |V(B)|  \end{equation}

Combining (\ref{eq:ht1}) and (\ref{eq:ht2}) we have that 
\begin{equation} \label{eq:ht3}
{\rm ht}(I_B)-{\rm ht}(I_{G'})= \begin{cases} 0,\ \textrm{if}\ G' \ \textrm{is bipartite} \\
1,\ \textrm{if}\ G' \ \textrm{is not bipartite} \end{cases} \end{equation}

Since every walk in $G'$ corresponds to a walk in $B$, we have that $I_{G'}\subseteq I_B$. Let us prove the converse statement. We know that $I_B$ is generated by quadrics, i.e., binomials which correspond to cycles of length four (see Corollary \ref{quadratic=allminimal=4cycles}) . Let $c=(v_1,v_2,v_3,v_4)$ be a cycle of $B$, we are going to build a cycle $c'$ in $G'$ such that $B_c = B_{c'}$. If $c$ does not pass through $x$, then we take $c' = c$. In case that $c$ passes through $x$, then it has the form $(x,v_1,v_2,v_3,x)$. Since the vertices of $c$ form a $K_{2,2}$, then there exists an $i \in \{1,\ldots,k\}$ such that $\{x,v_1,v_2,v_3,x\} \subseteq V(G_i)$. If $i = 1$ we choose  $c' := (x_1,v_1,v_2,v_3,x_1)$ and, if $i \geq 2$ we choose $c' := (x_2,v_1,v_2,v_3,x_2)$. It follows then that $I_{G'}\subseteq I_B$ and, therefore, \begin{equation} \label{eq:ht4} {\rm ht}(I_B)={\rm ht}(I_{G'}) \end{equation}

Combining (\ref{eq:ht3}) and (\ref{eq:ht4}) we conclude that $G'$ is bipartite. Since the ideal $I_G$ is generalized robust and generated by quadrics, by Corollary \ref{hereditary} it follows that $I_B=I_{G'}$ so is. Since $G'$ is bipartite, by Theorem \ref{quadratic-bipartite} all the blocks of $G'$ are of type $K_{2,\ell}$ or cut edges. By construction of the graph $G'$ we have that the block $B$ is a necklace graph of bipartite blocks each of them is of type $K_{2,\ell}$. Note also that every path from $x_1$ to $x_2$ has odd length because $B$ is not bipartite and $G$ is.  

\begin{figure}[h]
\includegraphics[scale=.6]{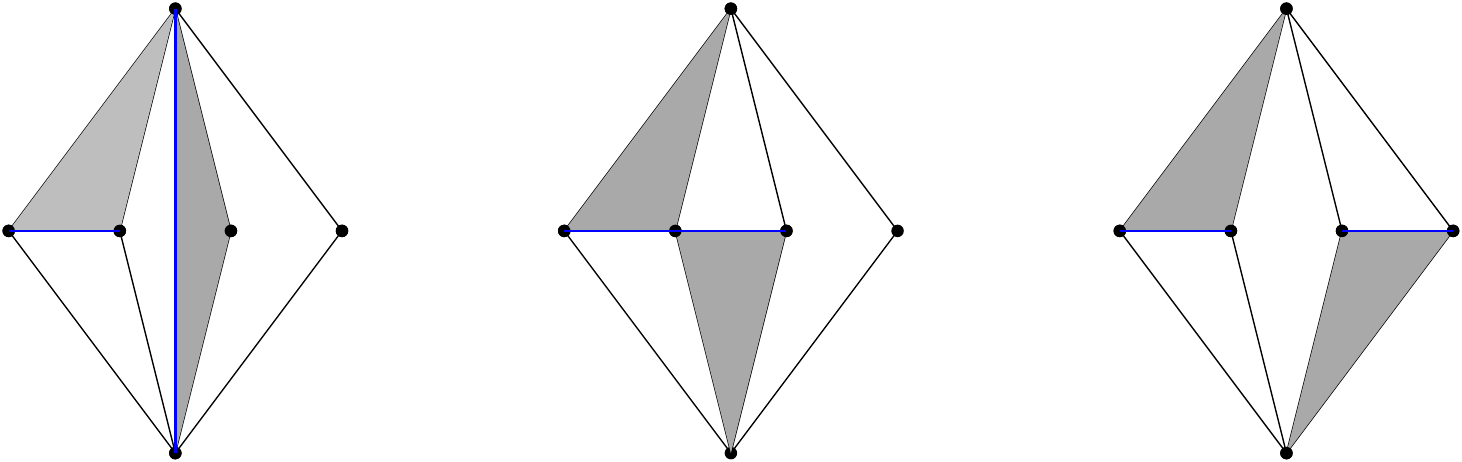}
\caption{The three non-isomorphic graphs $K_{2,4}$ plus two edges}
\label{fig:3cases}
\end{figure}

$(ii_b)$  We assume that there exists $i\in\{1,\ldots,k\}$ such that the graph $[G_i]$ is not bipartite. Then $[G_i]$ is a $K_{2,\ell}$ graph plus at least one more edge. We observe that we can only have one more edge or  $[G_i] = K_4$, otherwise we would be in one of the three cases shown in Figure \ref{fig:3cases}. In the three cases there are two edge disjoint triangles, which contradicts Lemma \ref{lm:shortcycles}.

\begin{figure}[h]
\includegraphics[scale=.6]{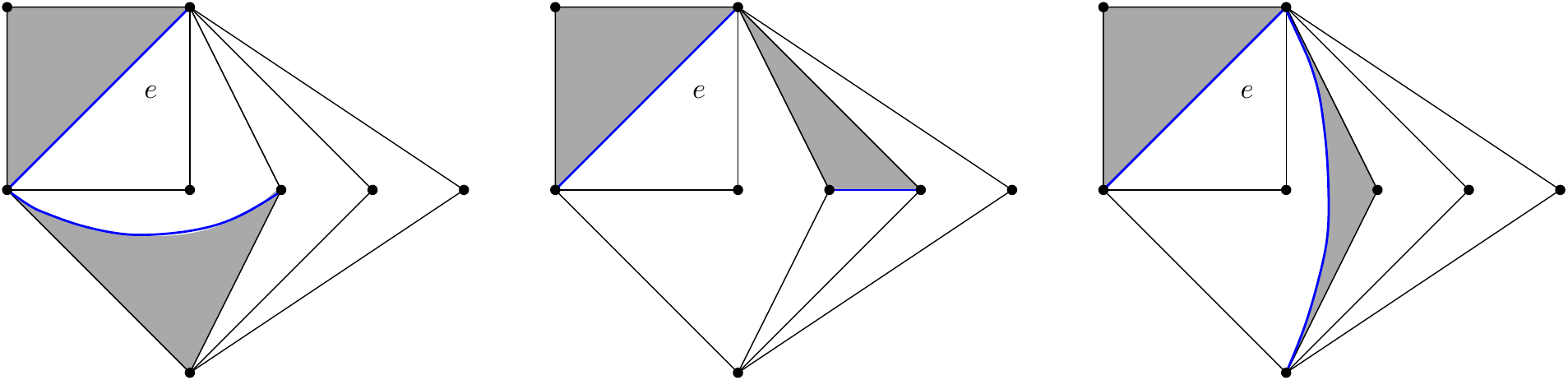}
\caption{The two $[G_i]$'s with a common edge are non bipartite}
\label{fig:3morecases}
\end{figure}

As a consequence, either i) $[G_i] = K_4,$ or ii) $[G_i] = K_{2,\ell} \cup \{e\}$. In ii), by (Claim $1$), there is another maximal $G_j$ such that $e \in E(G_j)$. We observe that $[G_j]$ is bipartite. Otherwise, there exists a subgraph as the one shown in Figure \ref{fig:3morecases}, again a contradiction to Lemma \ref{lm:shortcycles}.(b). 

Thus, in i) we have that $[G_i] = K_4$ and in ii) we have that the induced subgraph with vertex set $V(G_i) \cup V(G_j)$ is a double-$K_{2,(r,s)}$ graph. We just have to prove that in both cases the graph has no more vertices. We denote $H = K_4$ in i) and $H$ double-$K_{2,(r,s)}$ in ii). 
We observe that every two vertices in $H$ can be joined by two (non necessarily disjoint) paths: an odd path of length $\geq 3$ and an even path. Assume that there is another vertex $v$ at distance one of $u \in V(H)$. Since the graph is biconnected there exists a path from $v$ to a vertex of $V(H) - \{u\}$ and we take the shortest one. Hence, we have a path from $u$ to another vertex of $u' \in V(H)$ that only has its two endpoints in $V(H)$. If this path is odd, there is an even cycle of length $\geq 6$, a contradiction to Lemma \ref{lm:shortcycles}.  If this path is even, then there is even cycle of length $\geq 6$ or one of the $G_i$'s involved is not a maximal $K_{2,\ell}$ subgraph, a contradiction arises for both cases.

  ($\Longleftarrow$) We have exactly one non-bipartite block. As a consequence, every element in the Graver basis (and, hence, every minimal generator and every element in the universal Gr\"obner basis) corresponds to a walk entirely contained in a block. Thus, by Theorem \ref{quadratic-bipartite} it suffices to prove that $K_4$, $K_{2,\ell} \cup \{e\}$, a necklace-$K_{2,\ell}$ and a double-$K_{2(r,s)}$ graph give rise to generalized robust toric ideals generated by quadrics. Clearly $I_{K_4}$ is generalized robust  and generated by quadrics (see Example \ref{example}) and, since $I_{K_{2,\ell} \cup \{e\}} = I_{K_{2,\ell}}$ (they are both prime ideals, have the same height and one is contained in the other, so they are equal), so is $I_{K_{2,\ell} \cup \{e\}}$. The remaining two cases follow from Propositions \ref{pr:theguy} and \ref{pr:necklace}, and the proof is complete.
 \qed
 
 \

\begin{figure}[h]
\begin{center}
\includegraphics{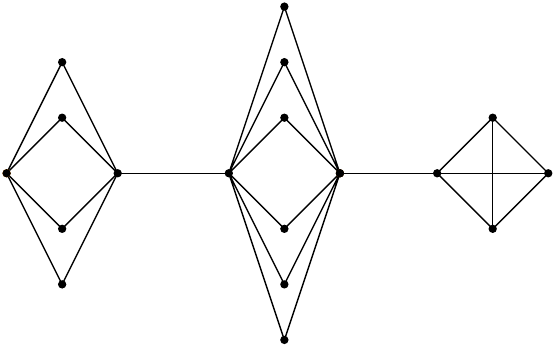}
\caption{A graph $G$ such that $I_G$ is quadratic generalized robust and not quadratic robust}
\label{genrobust}
\end{center}
\end{figure}

In Figure \ref{genrobust} we see an example of a graph $G$ whose corresponding toric ideal $I_G$ is a generalized robust ideal and is generated by quadrics. The graph $G$ consists of five blocks$;$ two cut edges, a $K_{2,4}$, a $K_{2,6}$ and exactly one non bipartite block which is a $K_4$. The existence of the $K_4$ as a subgraph, as we show in the next section, implies that the ideal $I_G$ is not robust.

\subsection{Quadratic robust graphs}
The goal of this section is to prove Theorem \ref{maintheorem2}, where we characterize the graphs providing robust toric ideals generated by quadrics. This result completes Corollary 5.2 of Boocher et al. in \cite{BOO2}. The main ingredients for proving the present theorem are Theorem \ref{maintheorem} and the following property:

\begin{thm1}\cite[Theorem 5.10]{CHT}\label{unique-robust} Let $I_A$ be a robust toric ideal. Then $I_A$ has a unique minimal system of generators.
\end{thm1}

From the above result it follows that a toric ideal $I_A$ with a unique minimal system of generators (i.e., $M_A=\mathcal{M}_A$) is robust if and only if it is a generalized robust (\cite[Corollary 5.13]{CHT}). 
From \cite{TT1} we know that $I_G$ is generated by indispensable binomials (i.e., it has a unique minimal system of generators)
 if and only if no closed walk $w$ such that $B_w$ is a minimal generator of $I_G$ has an $F_4$ (being the concept of $F_4$ rather technical, we refer to \cite{TT1} for its definition).
In the particular case that $I_G$ is generated by quadrics, the existence of an $F_4$ in a closed walk $w$ providing a minimal generator $B_w$ is equivalent to the existence of a subgraph $K_4$. Hence, the only obstruction for a toric ideal generated by quadrics to have a unique minimal set of generators is the existence of a subgraph $K_4$ (or, equivalently, having clique number $\geq 4$). We work out the example of $G = K_4$ in detail to show that $I_{K_4}$ is generalized robust but not robust.

\begin{figure}[h]
\begin{center}
\includegraphics{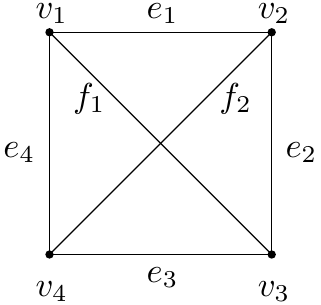}
\caption{The graph $K_4$ in Example \ref{example}.}
\label{fig:k4}
\end{center}
\end{figure}

\begin{ex1}\label{example}{\rm Consider the complete graph on four vertices $K_4$ on the vertex set $V(K_4)=\{v_1,v_2,v_3,v_4\}$ and on the edge set $E(K_4)=\{e_1=\{v_1,v_2\}, e_2=\{v_2,v_3\}, e_3=\{v_3,v_4\}, e_4=\{v_4,v_1\}, f_1=\{v_1,v_3\}, f_2=\{v_2,v_4\}\}$ (see Figure \ref{fig:k4}). It is clear that we have exactly three $4$-cycles $w_1=(v_1,v_2,v_3,v_4,v_1),\, w_2 = (v_1,v_2,v_4,v_3,v_1)$ and $w_3 = (v_1,v_4,v_2,v_3,v_1).$

The corresponding ideal $I_{K_4}$ is generated by the three binomials: $$I_{K_4}=\langle B_{w_1}=e_1e_3-e_2e_4, B_{w_2}=e_1e_3-f_1f_2, B_{w_3}=e_2e_4-f_1f_2\rangle.$$ Obviously, none of $B_{w_1}, B_{w_2}, B_{w_3}$ is indispensable since $$B_{w_i}\in \langle B_{w_j}, B_{w_k}\rangle, \ \textrm{for all distinct}\ i,j,k\ \textrm{where}\ i,j,k\in \{1,2,3\}.$$ Thus, the ideal has three different Markov bases$;$ 
 $$M_1=\{ B_{w_1}, B_{w_2}\}, M_2=\{ B_{w_1}, B_{w_3}\}, M_3=\{ B_{w_2}, B_{w_3}\}.$$ The universal Markov basis of the ideal is $\mathcal{M}_{K_4}=\{B_{w_1}, B_{w_2}, B_{w_3}\}$. It is easy to check that the universal Gr\"obner basis of the ideal $I_{K_4}$ is $\mathcal{U}_{K_4}=\{B_{w_1}, B_{w_2}, B_{w_3}\}$. It follows that the quadratic ideal $I_{K_4}$ is generalized robust but not robust. 
 }
 \end{ex1}

Putting all together we have that if $I_G$ is generated by quadrics, then $I_G$ is robust if and only if $I_G$ is generalized robust and does not have $K_4$ as subgraph. Thus, from Theorem \ref{maintheorem2} we deduce the following structural result.

\begin{thm1} \label{maintheorem2}Let $G$ be a non bipartite graph. The ideal $I_G$ is robust and is generated by quadrics if and only if all the blocks of $G$ are bipartite except one which is either a $K_{2,\ell}\cup \{e\}$ or a double-$K_{2,(r,s)}$ or a necklace-$K_{2,\ell}$ graph. The bipartite blocks of the graph $G$ are of type $K_{2,\ell}$ or cut edges.
\end{thm1}

\begin{rem1}{\rm In \cite{SUL} Sullivant introduces and studies the notion of strongly robust toric ideals. The motivation for studying strongly robust toric ideals comes from algebraic statistics. A toric ideal is strongly robust if its Graver basis coincides with its set of indispensable binomials, for more see \cite{SUL}. From \cite{CHT} one has that the notion of strongly robust and robust coincide for toric ideals of graphs. It follows that Theorem \ref{maintheorem2}, also characterizes completely the strongly robust quadratic toric ideals of graphs.
}
\end{rem1}

\begin{ex1}{\rm In Figure \ref{bothgenrobust-robust} we consider a graph $G$ which consists of six blocks$;$ two cut edges, a $K_{2,2}$, a $K_{2,4}$, a $K_{2,6}$ and one necklace-$K_{2,\ell}$ graph. By Theorem \ref{maintheorem} and Theorem \ref{maintheorem2} the corresponding ideal $I_G$ is both a robust and a generalized robust toric ideal and is generated by quadrics. The graph in Figure \ref{genrobust} contains a $K_4$ as a subgraph. The corresponding ideal is a quadratic generalized robust toric ideal, however it is not robust because of the existence of a $K_4$.

\begin{figure}[h]
\begin{center}
\includegraphics{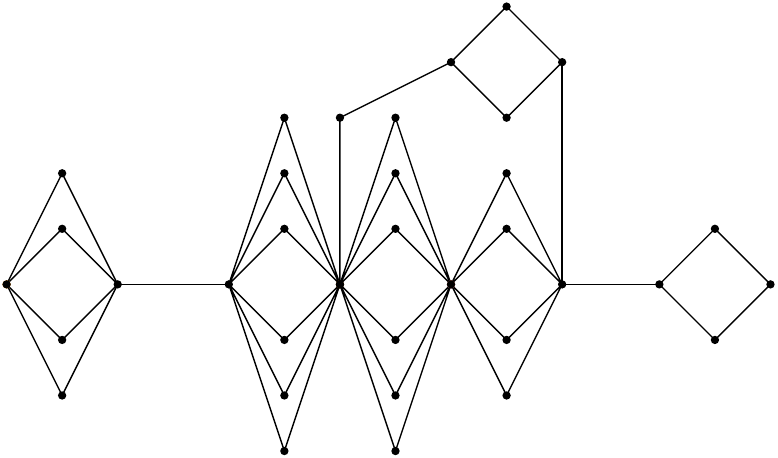}
\caption{A graph $G$ such that $I_G$ is both quadratic robust and quadratic generalized robust}
\label{bothgenrobust-robust}
\end{center}
\end{figure}
}
\end{ex1}

\section{Generalized robust ideals and numerical semigroups} \label{sec:numerical}

\subsection{Numerical semigroups having a complete intersection initial ideal}

A numerical semigroup is a submonoid $\mathcal S$ of $(\mathbb N,+)$ with finite complement. Every numerical semigroup has a unique minimal generating set  $A = \{a_1,\ldots,a_m\} \subseteq \Z^+$ of relatively prime integers, which is always finite. The number $m$ of elements of $A$ is usually called the {\it embedding dimension} of $\mS$.  The only numerical semigroup with embedding dimension $1$ is $\mS = \mathbb N$, here $A = \{1\}$ and $\IS$ is the zero ideal. From now on we assume that $\mS \subsetneq \N$ and, as a consequence, its embedding dimension is at least two.  
When we write $\mS = \langle a_1,\ldots,a_m \rangle$  we implicitly mean the numerical semigroup $\mS = \left\{ \sum_{i = 1}^m \alpha_i a_i \, \vert \, \alpha_i \in \N \right\}$ with minimal set of generators $A = \{a_1,\ldots,a_m\}$, and we call $\IS := I_A$ the toric ideal of the corresponding semigroup. 

By Krull's dimension theorem, any set of generators of an ideal $J \subseteq \mathbb{K}[\x]$ has at least ${\rm ht}(J)$ elements and is a complete intersection when equality occurs. Whenever $\mS = \langle a_1,\ldots,a_m \rangle$ is a numerical semigroup, then  $\IS$ has height $m-1$ and  $\IS$ is a complete intersection if and only if $\mu(\IS) = m-1$ or, in other words, if it can be generated by a set of $m-1$ binomials. 

Given a monomial ordering $\prec$, we have that ${\rm ht}({\rm in}_\prec(\IS)) = {\rm ht}(\IS)$ and that $\mu({\rm in}_\prec(\IS)) \geq \mu(\IS)$. Hence, whenever there exists a monomial ordering such that ${\rm in}_\prec(\IS)$ is a complete intersection, then so is $\IS$. 
In Theorem \ref{th:existsInitial} we characterize all numerical semigroups having a complete intersection initial ideal. 

A numerical semigroup $\mS$ with minimal generating set $A = \{a_1,\ldots,a_m\}$ is said to be {\it free for the arrangement} $a_1,\ldots,a_m$ if
 \begin{equation}\label{eq:free}  \lcm(a_i,\gcd(a_{i+1},\ldots,a_m)) \in \langle a_{i+1},\ldots, a_m\rangle, \text{ for all } i \in \{1,\ldots,m-1\}.
 \end{equation} Equivalently, if $\mu_i := \lcm(a_i,\gcd(a_{i+1},\ldots,a_m))$, the numerical semigroup $\mS$ is free for the arrangement $a_1,\ldots,a_m$ if there exist $\alpha_{(i,i+1)},\ldots,\alpha_{(i,m)}\in\mathbb{N}$ such that $\mu_i=\alpha_{(i,i+1)}a_{i+1}+\cdots+\alpha_{(i,m)}a_{m}$ for all $i \in \{1,\ldots,m-1\}$. We say that $\mS$ is {\it free} if it is free for an arrangement of its minimal generating set. 
 
 \begin{ex1}\label{ex:free}{\rm Consider the numerical semigroup $\mS = \langle a_1,a_2,a_3,a_4 \rangle$ with $a_1 = 8,\, a_2 = 9,\, a_3 = 10,\, a_4 = 12$. We have that $\mS$ is not free for the arrangement $a_1,a_2,a_3,a_4$ because $\lcm(a_1,\gcd(a_2,a_3,a_4)) = 8 \notin \langle a_2,a_3,a_4\rangle$. However, $\mS$ is free for the arrangement $a_2 = 9,\, a_3 = 10,\, a_1 = 8,\, a_4 = 12$. Indeed, \begin{itemize} 
 \item $\lcm(a_2,\gcd(a_1,a_3,a_4)) = 18 = a_1 + a_3 \in \langle a_1,a_3,a_4\rangle$,
 \item $\lcm(a_3,\gcd(a_1,a_4)) = 20 = a_1 + a_4 \in \langle a_1,a_4 \rangle$, and
 \item $\lcm(a_1,a_4) = 24 = 2 a_4 \in \langle a_4 \rangle$.
 \end{itemize}
 Thus, $\mS$ is a free numerical semigroup. 
 }
 \end{ex1}
 
 Equivalently, this notion can be inductively defined as follows: a numerical semigroup $\mS$ is free if either $\mS = \langle 1 \rangle = \mathbb N$ or there exists an arrangement $A = \{a_1,\ldots,a_m\}$ of its minimal generators such that  $d a_1  \in \langle a_{2},\ldots, a_m\rangle$ and the numerical semigroup $\mS' = \langle a_2/d, \ldots,a_m/d\rangle$ is free, where $d = \gcd(a_{2},\ldots,a_m)$. The following result can be found in several places in the literature, see, e.g., \cite[Lemma 2.1 and Proposition 2.3]{BG} or \cite[Lemma 3.2]{M}.

\begin{prop1}\label{pr:prefree} Let $\mS = \langle a_1,\ldots,a_m\rangle$ be a numerical semigroup and set $d := \gcd(a_2,\ldots,a_m)$. If $da_1 = \alpha_2 a_2 + \cdots + \alpha_m a_m$ with $\alpha_2,\ldots,\alpha_m \in \mathbb N$, then 
\[ \IS = I_{\mathcal S'} \cdot \mathbb K[x_1,\ldots,x_m] + \langle x_1^{d} - x_2^{\alpha_2} \cdots x_m^{\alpha_m} \rangle,\]
where $\mS' = \langle a_2/d,\ldots, a_m/d \rangle$ and $I_{\mathcal S'} \subseteq \mathbb K[x_2,\ldots,x_m]$.  
\end{prop1}

Applying inductively Proposition \ref{pr:prefree} we get the following result, which explains how to use condition (\ref{eq:free}) to construct a minimal set of generators of $\IS$ when $\mS$ is free.

\begin{prop1}\label{pr:free} Let $\mS = \langle a_1,\ldots,a_m\rangle$ be a free semigroup for the arrangement $a_1,\ldots,a_m$. Consider $\beta_i := \lcm(a_i,\gcd(a_{i+1},\ldots,a_m))/a_i$ and $\alpha_{(i,j)} \in \mathbb N$ so that  
\[ \lcm(a_i,\gcd(a_{i+1},\ldots,a_m)) =  \sum_{j = i+1}^{m} \alpha_{(i,j)} a_j \text{ for all } i \in \{1,\ldots,m-1\}. \] Then, $ \IS = \langle f_1,\ldots,f_{m-1}\rangle$, where $f_i = x_i^{\beta_i}  - \prod_{j = i+1}^m x_j^{\alpha_{(i,j)}}$ for all $i \in \{1,\ldots,m-1\} $.
\end{prop1}

As a direct consequence of Proposition \ref{pr:free}, every free semigroup has a set of generators consisting of $m-1$ binomials and, thus, it is a complete intersection. When $m = 2$ every numerical semigroup is free. For $m=3$, Herzog proved that $\mS$ is free if and only if $\IS$ is a complete interesection, see \cite{H}. For $m \geq 4$, there are complete intersection semigroups which  are not free. For example $\mS = \langle 10,14,15,21\rangle$ is not free for any arrangement of the generators and $\IS$ is a complete intersection; indeed, one can check that $\{ x_1^3 - x_3^2, x_2^3 - x_4^2, x_1^2x_3 - x_2x_4\}$ is a Markov basis for $\IS$ (and it is also the universal Markov basis of $\IS$). 

 In this section we will use the following general fact about toric ideals, which we write here only for numerical semigroups (see, e.g., \cite{ST} or \cite{V}). Let $\mS = \langle a_1,\ldots,a_m \rangle$ be a numerical semigroup and consider the morphism of groups $$\rho: \Z^m \longrightarrow \Z,\ \textrm{induced by} \ \rho(\mathbf{e}_i) = a_i,\ \forall i \in \{1,\ldots,m\},$$ being $\{\mathbf{e}_1,\ldots,\mathbf{e}_m\}$ the canonical basis of $\Z^m$.  Let $\mathbf{u}, \mathbf{v} \in \mathbb Z^m$ and consider the binomial $f = \mathbf x^{\mathbf{u}} - \mathbf x^{\mathbf{v}}$, we set $\widetilde{f} := \mathbf{u} - \mathbf{v} \in \Z^m$. We have that $f \in \IS$ if and only if $\widetilde{f} \in {\rm ker}(\rho)$. Moreover, we have the following:
\begin{prop1}\label{pr:rho} Let $\{f_1,\ldots,f_r\}$ be a set of binomials. If $\IS = \langle f_1,\ldots,f_r \rangle$, then ${\rm ker}(\rho) = \langle \widetilde{f}_1,\ldots, \widetilde{f}_r \rangle$. 
\end{prop1}

From the following two propositions easily follows the proof of Theorem \ref{th:existsInitial}, which is the main result of this subsection.

\begin{prop1}\label{pr:freeGB} Let $\mS$ be a numerical semigroup with minimal generating set $A = \{a_1,\ldots,a_m\}$. Then,
 $\mS$ is free for the arrangement $a_1,\ldots,a_m$ if and only if the reduced Gr\"obner basis with respect to the lexicographic order with $x_1 \succ \cdots \succ x_m$ has $m-1$ elements.
\end{prop1}
\begin{proof}$(\Longrightarrow)$ Suppose that  $\mS = \langle a_1,\ldots,a_m\rangle$ is free for the ordering $a_1,\ldots,a_m$. By Proposition \ref{pr:free} we have that  $$\IS = \langle f_1,\ldots,f_{m-1} \rangle,\ \textrm{where}\ f_i = x_i^{\beta_i}  - \prod_{j = i+1}^m x_j^{\alpha_{(i,j)}},\ \textrm{for some}\ \beta_i, \alpha_{(i,j)} \in \mathbb N.$$ Considering $\prec$ the lexicographic order with $x_1 \succ \cdots \succ x_m$ we observe that ${\rm in}_\prec(f_i)  = x_i^{\beta_i}$ for all $i \in \{1,\ldots,m-1\}$. Since the initial forms are pairwise prime, then $\mathcal G = \{f_1,\ldots,f_{m-1}\}$ is a Gr\"obner basis of $\IS$ for $\prec$ and, hence, the reduced Gr\"obner basis with respect to $\prec$ has $m-1$ elements. 

$(\Longleftarrow)$ Consider $\prec$ the lexicographic order with $x_1 \succ \cdots \succ x_m$ and let $\mathcal G$ be the corresponding reduced Gr\"obner basis of $\IS$. For all $i \in \{1,\ldots,m-1\}$ we have that $x_i^{a_{i+1}} - x_{i+1}^{a_i} \in \IS$ and, hence, $x_i^{a_{i+1}}  \in {\rm in}_\prec(\IS)$. As a consequence,   ${\rm in}_\prec(\IS) = \langle x_1^{b_1},\ldots,x_{m-1}^{b_{m-1}} \rangle$, with $b_1,\ldots,b_{m-1}\in\mathbb{Z}^+$ and $\mathcal G = \{g_1,\ldots,g_{m-1}\}$ with $g_i = x_i^{b_i} -  M_i$, being $M_i$ a monomial not involving the variables $x_1,\ldots,x_{i}$.

Our next goal is to prove that $b_1 = \gcd(a_{2},\ldots,a_m)$ from which we conclude that  $\lcm(a_1,\gcd(a_{2},\ldots,a_m)) \in \langle a_{2},\ldots, a_m\rangle$. 

Set $B := \gcd(a_2,\ldots,a_m)$, we observe that $$\Z B a_1 = \Z a_1 \cap (\sum_{j = 2}^m \Z a_j)$$ and, in particular, 
\begin{equation}\label{eq:B} B a_1 = \sum_{i = 2}^m \gamma_j a_j \text{ for some } \gamma_j \in \mathbb Z, \end{equation} Since $g_1 \in \IS$, it follows that $$b_1 a_1 = \deg_A(M_1) \in \Z a_1 \cap (\sum_{j = 2}^m \Z a_j)$$ and, then, $B$ divides $b_1$. 
Moreover, following the notation of Proposition \ref{pr:rho}, we have that
$${\rm ker}(\rho) = \langle \widetilde{g}_1,\ldots, \widetilde{g}_{m-1} \rangle\ \textrm{where}\ \widetilde{g}_i = b_i e_i -  \sum_{j > i}c_{(i,j)} e_j \in \Z^m$$
for some $c_{(i,j)}\in \mathbb N$ and, from (\ref{eq:B}), we deduce that $$B e_1 - \sum_{j = 2}^m \gamma_j e_j \in {\rm ker}(\rho).$$
Since $\widetilde{g}_1$ is the only element among $\widetilde{g_1},\ldots,\widetilde{g}_{m-1}$ with a nonzero first entry, we conclude that $b_1$ divides~$B$.

We have thus already proved that $$b_1 =  \gcd(a_2,\ldots,a_m)\ \textrm{and} \ g_1 = x_1^{b_1} - M_1 \in \IS.$$ Then, $$ a_1 \gcd(a_2,\ldots,a_m)  = \deg_A(M_1) \in \langle a_2,\ldots,a_m\rangle.$$ Since $\gcd(a_1,\gcd(a_2,\ldots,a_m)=1$ it follows that $$a_1 \gcd(a_2,\ldots,a_m) = \lcm(a_1,\gcd(a_{2},\ldots,a_m)) \in \langle a_{2},\ldots, a_m\rangle.$$

 Finally, it suffices to observe that $$\mathcal G \cap \mathbb K[x_2,\ldots,x_m] = \{g_2,\ldots,g_{m-1}\}$$ is the reduced Gr\"obner basis of $\IS \cap \mathbb K[x_2,\ldots,x_m] = I_{\mS'}$ with respect to the lexicographic order with $x_2 \succ \cdots \succ x_m$ and proceed inductively to get the  result.
\end{proof}

\begin{prop1}\label{pr:freeGB2} Let $\mS$ be a numerical semigroup with minimal generating set $A = \{a_1,\ldots,a_m\}$. If $\mathcal G = \{g_1,\ldots,g_{m-1}\}$ is a Gr\"obner basis of $\IS$ with respect to a monomial ordering $\prec$, then $\mathcal G$ is also a Gr\"obner basis with respect to a certain lexicographic monomial ordering $\prec_{lex}$.
\end{prop1}
\begin{proof}Since $\mathcal G$ has $m-1$ elements, then $\mu({\rm in}_\prec(\IS)) = m-1$. For all $i,j \in \{1,\ldots,m\}$ with $i\neq j$, we have that $x_i^{a_j} - x_j^{a_i} \in \IS$ and, hence, either $x_i^{a_j}$ or $x_j^{a_i}$ belongs to the initial ideal ${\rm in}_\prec(\IS)$. As a consequence,  we may assume (after reindexing the variables if necessary) that ${\rm in}_\prec(\IS) = \langle x_1^{b_1},\ldots,x_{m-1}^{b_{m-1}} \rangle$ for some $b_1,\ldots,b_{m-1}\in\mathbb{Z}^+$ and that $g_i = x_i^{b_i} -  M_i$, where $M_i$ is a monomial not involving the variable $x_i$, in which $i=1,\ldots,m-1$.

{\it Claim:} there exists an $\ell \in \{1,\ldots,m-1\}$ such
that $x_\ell$ does not divide $M_i$ for all $i \in \{1,\ldots,m-1\}$.

{\it Proof of the claim: } Assume by contradiction that the claim does not hold, i.e., for all $\ell\in \{1,\ldots,m-1\}$ there exists $i\in\{1,\ldots,m-1\}$ such that $x_\ell$ divides $M_i$. We consider the simple directed graph
with vertex set $\{1,\ldots,m-1\}$ and arc set $\{(j,i)\, \vert\, 1
\leq i,j \leq m-1$ and $x_j$ divides $M_i\}$. Then, the out-degree of every vertex is greater or equal to one and, thus, there is a directed cycle in the graph. Assume, without
loss of generality, that the cycle is $(1,2,\ldots,r,1)$
with $r \leq m-1$. This implies that $g_i\in\langle x_{i-1},x_i\rangle,\ \forall i=1,\ldots,r-1$ and $x_0=x_r$, thus: $$\langle g_{1},\ldots,g_{r}\rangle \subsetneq
\langle x_{1},\ldots,x_{r}\rangle,$$ and so $$\IS \subsetneq H :=
\langle x_{1},\ldots,x_{r},g_{r+1},\ldots,g_{m-1}\rangle,$$ but this is not
possible because $\IS$ is prime and $$m-1 = {\rm ht}(\IS) < {\rm
ht}(H) \leq m-1.$$

Hence, $x_\ell$ only appears in $g_\ell$. Assume without loss of generality that $\ell = 1$.   
Proceeding as before, one can prove that there exists an $\ell' \in \{2,\ldots,m\}$ such that $x_\ell'$ does not divide $M_i$ for all $i \in \{2,\ldots,m-1\}$. Iterating this idea one gets that, after reindexing the variables if necessary, the variables $x_1,\ldots,x_i$ do not divide $M_i$. Hence, taking $\prec_{lex}$ the lexicographic order with $x_1 \succ_{lex} \cdots \succ_{lex} x_m$ one has that ${\rm in}_{\prec_{lex}}(g_i) = x_i^{b_i}$. Since they are all relatively prime, then $\mathcal G$ is also a Gr\"obner basis for $\prec_{lex}$. 
\end{proof}

Now, we can prove the main result of this subsection.

 \begin{thm1}\label{th:existsInitial} Let $\mS$ be a numerical semigroup. Then,
 $\mS$ is free if and only if it has a Gr\"obner basis with $m-1$ elements.
\end{thm1}
\begin{proof}$(\Longrightarrow)$ Follows directly from Proposition \ref{pr:freeGB}.

$(\Longleftarrow)$ Assume that $\mathcal G = \{g_1,\ldots,g_{m-1}\}$ is a Gr\"obner basis of $\IS$ with respect to a monomial ordering $\prec$. By Proposition \ref{pr:freeGB2}, $\mathcal G$ is also a Gr\"obner basis with respect to a certain lexicographic monomial ordering $\prec_{lex}$. The result follows from Proposition \ref{pr:freeGB}. 
\end{proof}

Let us illustrate these result with some examples.

\begin{ex1}{\rm Consider the numerical semigroup $\mS = \langle a_1, a_2, a_3, a_4  \rangle$ with $a_1 = 8,\, a_2 = 9,\, a_3 = 10,\, a_4 = 12$ of Example \ref{ex:free}. Since $\mS$ is not free for the arrangement $a_1,a_2,a_3,a_4$, Proposition \ref{pr:freeGB} assures that the reduced Gr\"obner basis with respect to the lexicographic order with $x_1 \succ x_2 \succ x_3 \succ x_4$ has more than $3$ elements. Indeed, it has $8$ elements. Nevertheless, $\mS$ is free for the arrangement $a_2,a_3,a_1,a_4$ and, again by Proposition \ref{pr:freeGB}, we know that the reduced Gr\"obner basis with respect to the lexicographic order with $x_2 \succ x_3 \succ x_1 \succ x_4$ has $3$ elements. Indeed, it is $\{x_2^2 - x_1 x_3,\, x_3^2 - x_1 x_4, x_1^3 - x_4^2\}$, which also is a Markov basis of $\IS$ (and it is also the universal Markov basis of $\IS$). 
}
\end{ex1}

\begin{ex1}{\rm Consider the numerical semigroup $\mS = \langle a_1, a_2, a_3, a_4  \rangle$ with $a_1 = 10,\, a_2 = 14,\, a_3 = 15,\, a_4 = 21$. We know that $\mS$ is not free and $\IS$ is a complete intersection. Thus, by Theorem \ref{th:existsInitial}, we conclude that $\IS$ cannot be minimally generated by a Gr\"obner basis. 
}
\end{ex1}

\subsection{Numerical semigroups defining robust - generalized robust toric ideals}
  
In \cite{GOR}, Garc\'ia-S\'anchez, Ojeda and Rosales studied a family of affine submonoids of $\N^n$ which they called {\it semigroups with a unique Betti element}. 

\begin{def1} An affine monoid $\mS$ with minimal generating set $A$ has a unique Betti element if and only if the $A$-degrees of all the binomials in a Markov basis of $I_A$ coincide. 
\end{def1}

In \cite[Theorem 6]{GOR}, the authors characterized semigroups with a unique Betti element as those where $\mathcal C_A = Gr_A$, i.e., the set of  circuits of $I_A$ coincides with the Graver basis. Moreover, in this family of semigroups one has that every circuit is a minimal generator of $I_A$. As a consequence of these two facts and Proposition \ref{pr:inclusions}, one has that all the toric bases coincide ($\mathcal C_A = \mathcal M_A = \mathcal U_A = Gr_A$) and one directly derives the following result.
 
 \begin{prop1}\label{pr:uniquebetti} Every affine monoid with a unique Betti element defines a generalized robust toric ideal. 
 \end{prop1}
 
 The converse of this result does not hold for general toric ideals (in Section \ref{sec:graph} one can find examples of generalized robust toric ideals of graphs not having a unique Betti element). Nevertheless, in the following result, which is the main one in this section, we aim at proving that the converse of Proposition \ref{pr:uniquebetti} holds for numerical semigroups.

\begin{thm1}\label{th:genrbustnumsemigroup} A numerical semigroup defines a generalized robust toric ideal if and only if it has a unique Betti element.
\end{thm1}

In the proof of this result, we handle with Betti divisible numerical semigroups. This is a family of numerical semigroups studied in \cite{GH} that contains those with a unique Betti element.

\begin{def1}
A numerical semigroup $\mS$ (and more generally an affine monoid) with minimal generating set $A$ is {\it Betti divisible} if  
the $A$-degrees of all the binomials in a Markov basis of $I_A$ are ordered by divisibility.
\end{def1}

Our strategy for proving Theorem \ref{th:genrbustnumsemigroup} is the following. We first study the case $m = 3$ and prove that whenever $\mS = \langle a_1,a_2,a_3 \rangle$ satisfies that  $\mathcal C_{\IS} \subseteq \mathcal M_{\IS}$, then it is Betti divisible (Proposition  \ref{pr:casem3}).  Then, we move on to the case of a numerical semigroup $\mS = \langle a_1,\ldots,a_m\rangle$ defining a generalized robust toric ideal. Since $\IS$ is generalized robust (i.e. $\mathcal M_{\IS}=\mathcal U_{\IS}$), by Theorem \ref{pr:inclusions} we have that $\mathcal C_{\IS} \subseteq \mathcal M_{\IS}$. By Proposition \ref{genrobustfewvariables}, we have that for all $A' \subseteq \{a_1,\ldots,a_m\}$, if we take $\mS' = \langle A' \rangle$, then  $\mathcal C_{I_{S'}} \subseteq \mathcal M_{I_{S'}}$. In particular, if we take $A'$ a set of three elements, by Proposition \ref{pr:casem3}, we have that $\mS'$ is Betti divisible. By conveniently choosing the set $A'$ and using the fact that $\mathcal U_{\IS} = \mathcal M_{\IS}$ we will conclude that $\mS$ has a unique Betti element. 

We introduce some concepts and results that we will use in the proof. Firstly, we have that the set of the circuits of the toric ideal of a numerical semigroup is given by the following result (see, e.g., \cite[Chapter 4]{ST} or \cite[Lemma 2.8]{KO}).
\begin{lem1}\label{lm:circuitsnumerical}Let $\mS = \langle a_1,\ldots,a_m \rangle \subseteq \mathbb N$ be a numerical semigroup. Then,
\[ \mathcal C_{\IS} = \left\{ q_{i,j} := x_i^{a_j/\gcd(a_i,a_j)}  -  x_j^{a_i/\gcd(a_i,a_j)} \, \vert \, 1 \leq i, j \leq m,\, i \neq j  \right\}, \]

\end{lem1}

In the forthcoming we will use the concept of critical binomial, which was introduced by Eliahou \cite{E} and later studied in \cite{AV} and \cite{KO}, among others. Let $\mS = \langle a_1,\ldots,a_m\rangle$ be a numerical semigroup, one sets $$n_i ={\rm min}\left\{\ b \in \mathbb Z^+ \ \vert \ b a_i \in \sum_{j \in \{1,\ldots,m\} \setminus \{i\}} \mathbb N a_j\ \right\},\ \textrm{for}\ i=1,\ldots,m.$$
Write $$n_i a_i = \sum_{j \in \{1,\ldots,m\} \setminus \{i\}}  \beta_j a_j,\ \textrm{with}\ \beta_j \in \N,$$ the binomials  \[ g_i := x_i^{n_i} - \prod_{j \in \{1,\ldots,m\} \setminus \{i\}}  x_j^{\beta_j}\ \textrm{and}\ -g_i \] of $A$-degree $n_i a_i$ are called \emph{critical binomials with respect to $x_i$}.

As we mentioned before, a key point in the proof is the case of embedding dimension $m = 3$.  Three-generated numerical semigroups and their  toric ideals have been extensively studied in the literature. Here we will recall some results concerning them that we will use later; one can find restatements of these results and their proofs in \cite{AG, H}. 
Let $\mS = \langle a_1,a_2,a_3\rangle$ be a numerical semigroup, then $2 \leq \mu(\IS) \leq 3$ and   the $A$-degrees of $\IS$ are $\{n_1 a_1, n_2 a_2, n_3a_3\}$. Moreover, $\mu(\IS)  = 2$ or, equivalently, $\IS$ is a complete intersection if and only if there exist $1 \leq i < j \leq 3$ such  that $n_i a_i = n_j a_j$.  Clearly $\mS$ has a unique Betti element if and only if  $n_1 a_1 = n_2 a_2 = n_3 a_3$.

\begin{prop1}\label{pr:casem3} Let $\mS = \langle a_1,a_2,a_3 \rangle$ be a numerical semigroup. If $\mathcal C_{\IS} \subseteq \mathcal M_{\IS}$, then $\mS$ is Betti divisible.
\end{prop1}
\begin{proof} Suppose that $n_1 a_1 \leq n_2 a_2 \leq n_3 a_3$. 
We assume that  $\mathcal C_{\IS} \subseteq \mathcal M_{\IS}$ and let us see that $n_1 a_1 = n_2 a_2$ and that they both divide $n_3 a_3$. 

{\it Claim 1:} $n_1 a_1 = n_2 a_2$.

{\it Proof of claim 1:} if $n_1 a_1 < n_2 a_2$, we write $n_1 a_1 = \alpha_2 a_2 + \alpha_3 a_3$ with $\alpha_2,\alpha_3 \in \N$. Since $n_1 a_1 < n_2 a_2 \leq n_3 a_3$ it follows that both $\alpha_2$ and $\alpha_3$ are nonzero. Take the critical binomial $f := x_1^{n_1} -  x_2^{\alpha_2} x_3^{\alpha_3} \in \IS$. Consider now the circuit $$q_{1,2} = x_1^{a_2/\gcd(a_1,a_2)} - x_2^{a_1/\gcd(a_1,a_2)} \in \mathcal C_{\IS} \subseteq \mathcal M_{\IS}.$$ We have  that $a_2/\gcd(a_1,a_2) > n_1$ (otherwise, $a_2/\gcd(a_1,a_2) = n_1$, then $n_1 a_1$ is a multiple of $a_2$ and $n_1 a_1 \geq n_2 a_2$, a contradiction). Hence, $$q_{1,2} - x_1^{a_2/\gcd(a_1,a_2) - m_1} f = x_2 h$$ for some $h \in \IS$. Thus, $q_{1,2} \in \langle x_1,\ldots,x_m\rangle\cdot \IS$ and, by Proposition \ref{pr:minimal}, $q_{1,2} \notin \mathcal M_{\IS}$, a contradiction.

{\it Claim 2:} $n_3 a_3$ is a multiple of both $a_1$ and $a_2$.

{\it Proof of claim 2:}  Assume $n_3a_3$ is not a multiple of $a_1$, then $n_1 a_1 = n_2 a_2 < n_3 a_3 < \lcm(a_1,a_3)$. Take a critical binomial $f := x_3^{n_3} - x_1^{\gamma_1} x_2^{\gamma_2} \in \IS$ with respect to $x_3$. We observe that $n_3 < a_1/\gcd(a_1,a_3)$. Also, we may assume that $\gamma_1 > 0$, otherwise we have that $\gamma_2 > n_2$ and we consider $f = x_3^{n_3} - x_1^{n_1} x_2^{\gamma_2 - n_2} \in \IS$.  As a consequence, the circuit $$q_{3,1} = x_3^{a_1/\gcd(a_1,a_3)} - x_1^{a_3/\gcd(a_1,a_3)} = x_3^{a_1/\gcd(a_1,a_3)- n_3} f + x_1 h,$$ for some $h \in \IS$.  Thus, $q_{3,1} \in \langle x_1,\ldots,x_m\rangle\cdot \IS$ and, by Proposition \ref{pr:minimal}, $q_{3,1} \notin \mathcal M_{\IS}$, a contradiction.  Hence, $n_3 a_3$ is a multiple of $a_1$. A similar argument proves that $n_3 a_3$ is also a multiple of $a_2$.

Now, by (Claim 1), we have that $n_1 a_1 = n_2 a_2 = \lcm(a_1,a_2)$ and, by (Claim 2), $n_3 a_3$ is a multiple of $\lcm(a_1,a_2)$. Thus, $\mS$ is Betti divisible and the result holds.

\end{proof}

Now we can proceed to the proof of the main theorem of this subsection. 

{\it \noindent Proof of Theorem \ref{th:genrbustnumsemigroup}} $(\Longleftarrow)$ This is a particular case of Proposition \ref{pr:uniquebetti}.
 
$(\Longrightarrow)$ Let $\mS = \langle a_1,\ldots,a_m \rangle$ be a numerical semigroup defining a generalized robust toric ideal. Let us prove that it has a unique Betti element. Assume that $n_1 a_1 \leq n_2 a_2 \leq \cdots \leq n_m a_m$.

{\it Claim 1:} $n_1 a_1 = n_2 a_2$.

Let us first see that $n_1 a_1$ is a multiple of $a_i$ for some $i \in \{2,\ldots,m\}$. Assume, by contradiction that this is not true. We write 
$$n_1 a_1 = \sum_{k = 2}^m \alpha_k a_k\ \textrm{with}\ \alpha_2,\ldots,\alpha_m \in \mathbb N,$$ consider $$f := x_1^{n_1} - \prod_{j = 2}^m x_j^{\alpha_j} \in \IS$$  and observe that at least two of the $\alpha_i$'s are nonzero. We take $s \in \{2,\ldots,m\}$ such that $\alpha_s \neq 0$, and we are going to see that the circuit $$q_{s,1} = x_1^{a_s/\gcd(a_1,a_s)} - x_s^{a_1/\gcd(a_1,a_s)}$$ is not a minimal binomial, which contradicts the hypothesis.  We have that $n_1 < a_s/\gcd(a_1,a_s)$. We  set
\[h := q_{s,1} - x_1^{\frac{a_s}{\gcd(a_1,a_s)} - n_1} f \in \IS,\]
and we have that \[ h = x_1^{\frac{a_s}{\gcd(a_1,a_s)} - n_1} \prod_{j = 2}^m x_j^{\alpha_j} - x_s^{\frac{a_1}{\gcd(a_1,a_s)}} \neq 0\] and the two monomials apprearing in $h$ are multiples of $x_s$. Then, we have that $h = x_s h'$ for some $h' \in \IS$. We conclude that $q_{s,1} \in \langle x_1,\ldots,x_m\rangle \cdot \IS$ and, by Proposition \ref{pr:minimal}, $q_{s,1} \notin \mathcal M_{\IS}$, a contradiction.

So far we have seen that $n_1 a_1$ is a multiple of $a_i$ for some $i \in \{2,\ldots,m\}$. In particular, we have that $n_1 a_1 = \lcm(a_1,a_i) \geq n_i a_i$. Thus, $n_1 a_1 = n_2 a_2 = n_i a_i$ and the claim follows.

{\it Claim 2:} $n_1 a_1 = n_k a_k$ for all $k \in \{3,\ldots,m\}$. 

Take $k \in \{3,\ldots,m\}$ and let us prove that $n_1 a_1 = n_k a_k$. Assume by contradiction that $n_1 a_1 < n_k a_k$.  Set $n_k' = \min \{b \in \mathbb Z^+ \, \vert \, b a_k \in \langle a_1,a_2 \rangle \}$, we have that $n_1 a_1 < n_k a_k \leq n_k'a_k$. We consider the semigroup $\mS' = \langle a_1,a_2,a_k \rangle$. Since $\IS$ is generalized robust, then $\mathcal U_{\IS} \subseteq \mathcal M_{\IS}$ and, by Proposition \ref{genrobustfewvariables}.(a), we have $\mathcal U_{I_{\mS'}} \subseteq \mathcal M_{I_{\mS'}}$. Now, applying Proposition \ref{pr:casem3} we get that $\mS'$ is Betti divisible. Since the Betti elements of $\mS'$ are $n_1 a_1 = n_2 a_2 < n_k' a_k$ we have that $n_1 a_1 \mid n_k'a_k$ and, as a consequence, $n_k' a_k = \lcm(a_1,a_k) = b\, n_1 a_1 = b\, n_2 a_2$ for some $b \geq 2$. 

Now we have that the circuit \[ q_{1,k} = x_1^{a_k/\gcd(a_1,a_k)} - x_k^{a_1/\gcd(a_1,a_k)} = x_1^{b\, n_1} - x_k^{n_k'} \in \mathcal C_{\IS} \subseteq \mathcal U_{\IS} = \mathcal M_{\IS}.\]  Set $p := x_1^{n_1} x_2^{(b-1) n_2} - x_k^{n_k'} \in \IS$, the equality
\[ p = q_{1,k} - x_1^{(b-1)n_1} (x_1^{n_1} - x_2^{n_2}) \]
implies that $p \in \mathcal M_{\IS}$ (see Proposition \ref{pr:minimal}). Nevertheless, $x_1^{n_1} - x_2^{n_2} \in \IS$ and, hence, for every monomial order $\prec$, either $x_1^{n_1}$ or $x_2^{n_2}$ belongs to ${\rm in}_{\prec}(\IS)$. In both cases the monomial $x_1^{n_1} x_2^{(b-1) n_2} \in {\rm in}_{\prec}(\IS)$ but is not a minimal generator of ${\rm in}_{\prec}(\IS)$. Thus, $p \notin \mathcal U_{\IS}$, a contradiction. 

From (Claim 1) and (Claim 2) we conclude that $\mS$ has a unique Betti element.

\qed

The shape of the generators of a numerical semigroup with a unique Betti element is given by the following result (see \cite{KO} or \cite[Example 12]{GOR}):

\begin{prop1}\label{pr:shapeuniquqBetti} A numerical semigroup $\mS = \langle a_1,\ldots,a_m \rangle \subseteq \N$ has a unique Betti element if and only if  there exist $d_1,\ldots,d_m \geq 2$ pairwise prime integers such that $a_i = (\prod_{j =1}^m d_j)/d_i$ for all $i \in \{1,\ldots,m\}$.
\end{prop1}

Hence, one can easily construct examples of generalized robust toric ideals of numerical semigroups with arbitrarily large embedding dimension.
Since numerical semigroups with a unique Betti element only have a unique minimal set of generators when $m = 2$, we directly get the following.

 \begin{cor1}\label{cor:robust} A numerical semigroup $\mS$ defines a robust toric ideal if and only if  $\mS$ is $2$-generated.
 \end{cor1}

\section{Conclusions and open questions}
\label{sec:conclusion}

In this paper we have described different families of robust and generalized robust toric ideals. In the context of toric ideals of graphs generated by quadrics, the family of generalized robust toric ideals is only slightly bigger than the family of robust ideals. Nevertheless, for toric ideals of numerical semigroups the situation changes: while only principal ideals are robust, there are generalized robust ideals with arbitrarily large number of generators. 

In \cite{CHT}, the author asks if $\mathcal M_{A} = Gr_A$ for generalized robust toric ideals and verifies it for toric ideals of graphs. As a consequence of Theorem \ref{th:genrbustnumsemigroup}, we provide an affirmative answer in the case of numerical semigroups. For toric ideals of graphs the equality $\mathcal M_{A} = Gr_A$ characterizes generalized robustness (Theorem \ref{M=Gr}). However, this is not true for numerical semigroups. Indeed, consider the (Betti divisible) numerical semigroup $\mS = \langle a_1, a_2, a_3 \rangle$ with $a_1 = 10 = 2 \cdot 5,\, a_2 = 12 = 2^2 \cdot 3,\, a_3 = 15 = 3 \cdot 5$. By Proposition \ref{pr:shapeuniquqBetti} we have that $\mS$ is not a semigroup with a unique Betti element and, thus, it is not generalized robust (Theorem \ref{th:genrbustnumsemigroup}). Nevertheless, the equality $\mathcal M_{A} = Gr_A = \{x_1^3 - x_3^2, x_2^5 - x_1^6, x_2^5 - x_1^3x_3^2, x_2^5 - x_3^4\}$ holds. The same example also shows that the containment $\mathcal M_A \subseteq \mathcal U_A$, which holds for toric ideals of graphs, does not always work for toric ideals of numerical semigroups.

By \cite[Proposition 2.5]{BOO1}, we have that robustness property is preserved under an elimination of variables. However, we do not know if the same result is true when we replace robustness by generalized robustness. By Proposition \ref{genrobustfewvariables}.(a) we know that whenever $I_A$ is generalized robust and $A' \subseteq A$, then $\mathcal U_{A'} \subseteq \mathcal M_{A'}$, but we do not know if equality holds. 
 
\begin{quest1}\label{question}{\rm Let $I_A$ be a generalized robust toric ideal and $A' \subseteq A$, is $I_{A'}$ generalized robust?
}
\end{quest1}

In the present paper we give a positive answer to this question for: (1) toric ideals of graphs (see Corollary \ref{hereditary}), and (2) toric ideals of numerical semigroups (this follows as a consequence of Theorem \ref{th:genrbustnumsemigroup} and Proposition \ref{pr:shapeuniquqBetti}). 

In Theorem \ref{th:existsInitial} we characterize when the toric ideal of a numerical semigroup has a complete intersection initial ideal. It would be interesting to seek the answer to the same question for toric ideals of graphs. Since having a complete intersection initial ideal implies that the ideal itself is a complete intersection, the class of graphs that we are looking for is a subfamily of the one described in \cite{TT3}.

\begin{prob1}{\rm Characterize when the toric ideal $I_G$ of a graph $G$ has a complete intersection initial ideal.
}
\end{prob1}

Also, in Theorem \ref{th:existsInitial}, we proved that free numerical semigroups have an initial ideal such that $\mu(\IS) = \mu({\rm in}_{\prec}(\IS))$. There are further families of numerical semigroups with the same property. For example when $\mS$ is generated by an arithmetic sequence of integers (see, e.g., \cite{GSS}).

\begin{prob1}{\rm Characterize the numerical semigroups such that  $$\mu(\IS) = \mu({\rm in}_{\prec}(\IS))$$ for a  monomial order $\prec$. 
}
\end{prob1}

We have verified when the equality $\mathcal M_{\IS} = \mathcal U_{\IS}$ occurs for a numerical semigroup $\mS$. It would be interesting to characterize when equality or containment of other toric bases holds. For example, we say that a toric ideal $I_A$ is a {\it circuit ideal} if it is generated by its set of circuits (see \cite{BJT,MBV} for a deeper study of circuit ideals).  

\begin{prob1}{\rm Characterize the numerical semigroups $\mS$ such that $\IS$ is a circuit ideal.
}
\end{prob1}

We do not know the answer to this question even if we add the hypothesis of $\mS$ being a complete intersection. Indeed, whenever $\mS$ is a Betti divisible numerical semigroup, then it is a complete intersection and $\IS$ is a circuit ideal  (see \cite[Section 7]{GH}). However, there are further examples of complete intersection numerical semigroups such that $\IS$ is generated by $\mathcal C_{\IS}$. For example, consider the 
numerical semigroup $\mS = \langle a_1,a_2,a_3,a_4 \rangle$ with $a_1 = 390,\, a_2 =  546,\, a_3 =  770,\, a_4 = 1155\rangle$. Then, $\IS$ is minimally generated by $\{x_1^7 - x_2^5, x_3^3 - x_4^2, x_2^{55} - x_4^{26}\}$ and, thus, it is a complete intersection and a circuit ideal. However, the Betti degrees of the generators are $\beta_1 = 7 a_1 = 5a_2 = 2730,\, \beta_2 = 3 a_3 = 2 a_4 = 2310$ and $\beta_3 = 55 a_2 = 26 a_4 = 30030$ and, hence, it is not Betti divisible. Interestingly, in the context of toric ideals of graphs, every complete intersection ideal is a circuit ideal (by \cite[Theorem~5.1]{TTci}). 

\medskip

\noindent {\bf Acknowledgments}
\newline
This paper was written during the visit of the second author at the Department of Mathematics in Universidad de La Laguna (ULL). 
This work was partially supported by the Spanish MICINN ALCOIN (PID2019-104844GB-I00) and by the ULL funded research projects MASCA and MACACO.

Computational experiments with the computer softwares CoCoA \cite{ABR} and Singular \cite{SINGULAR} have helped in the elaboration of this work. 

\medskip
\noindent {\bf{Data Availability Statements}}
\newline
All data generated or analysed during this study are included in this published article (and its supplementary information files).

\end{document}